\documentclass[11pt,twoside]{article}
\usepackage{amssymb,amsmath,amsthm,mathrsfs,txfonts, color}
\newtheoremstyle{mythm}{3pt}{3pt}{}{}{\bfseries}{}{5mm}{}
\theoremstyle{mythm}
\def\pd#1#2{\frac{\partial#1}{\partial#2}}

\newtheorem{theorem}{\indent Theorem}[section]
\newtheorem{lemma}{\indent Lemma}[section]

\newtheorem{definition}{\indent Definition}[section]
\newtheorem{remark}{\indent Remark}[section]

\allowdisplaybreaks
\begin{document}

\title{Early and late stage profiles for a new chemotaxis model with
density-dependent jump probability and quorum-sensing mechanisms}

\author{
Tianyuan Xu$^a$, Shanming Ji$^b$\thanks{Corresponding author, email:jism@scut.edu.cn}, Chunhua Jin$^a$, Ming Mei$^{c,d}$,  Jingxue Yin$^a$,
\\
\\
{ \small \it $^a$School of Mathematical Sciences, South China Normal University}
\\
{ \small \it Guangzhou, Guangdong, 510631, P.~R.~China}
\\
{ \small \it $^b$School of Mathematics, South China University of Technology}
\\
{ \small \it Guangzhou, Guangdong, 510641, P.~R.~China}
\\
{ \small \it $^c$Department of Mathematics, Champlain College Saint-Lambert}
\\
{ \small \it Quebec,  J4P 3P2, Canada, and}
\\
{ \small \it $^d$Department of Mathematics and Statistics, McGill University}
\\
{ \small \it Montreal, Quebec,   H3A 2K6, Canada}
}
\date{}
\maketitle

\begin{abstract}
In this paper,
we derive a new chemotaxis model with degenerate diffusion
and density-dependent chemotactic sensitivity,
and we provide a more realistic description of cell migration process
for its early and late stages.
Different from the existing
studies focusing on the case of non-degenerate diffusion, the new model with
degenerate diffusion causes us some essential difficulty on the boundedness estimates
and the propagation behavior of its compact support.
In the presence of logistic damping, for the early stage before tumour cells spread to the whole body,
we first  estimate  the expanding speed of tumour region as $O(t^{\beta})$ for $0<\beta<\frac{1}{2}$.
Then, for the late stage of cell migration, we further prove that the asymptotic profile of the original system is just its corresponding steady state.
The global convergence of the original weak solution to the steady state with exponential rate $O(e^{-ct})$ for some $c>0$ is also obtained.
\end{abstract}

\baselineskip=18pt

\section{Introduction}
The motion of cells moving towards the higher concentration
of a chemical signal is called chemotaxis.
For example, bacteria moves toward the highest
concentration of food molecules to find food.
A well-known chemotaxis model was initially
proposed by Keller and Segel \cite{Keller} in 1971, subsequently,
a number of variations of the Keller-Segel system were proposed
and have been extensively studied  during the past four decades, for example,
see the survey papers \cite{Bellomo,Hillen} and
the references therein.
Especially, chemotaxis models also appear in medical mathematics.
Many factors effect the migration mechanisms
of tumour cells.
For example,
the extracellular matrix (ECM),
to which the tumour cell to be attached,
inhibits the cell polarizes and elongates to migrate.
ECM-degrading enzymes (MDE) cleave ECM fibers
into smaller chemotactic fragments to facilitate cell-migration \cite{Friedl}.
In \cite{Chaplain},
Chaplain and Anderson
introduced a model for tumour invasion mechanism,
which describes tumour invasion phenomenon
in accounting for the role of chemotactic ECM fragments named ECM*,
produced by a biological reaction
between ECM and MDE.
In these models, the cancer cell random motility is assumed
to be a constant,
which leads to linear isotropic diffusion.
However, in realistic situation, it is emphasized in that migration
of the cancer cells through the ECM fibers should rather be regarded like movement
in a porous medium with degenerate diffusion from a physical point of view \cite{Tao}.
Compared with the classical tumour invasion model with linear diffusion,
the mathematical analysis of the nonlinear diffusion system has to cope
with considerable additional challenges and is much less understood.
Several chemotaxis models with nonlinear  diffusion have been recently
proposed and analyzed, e.g. \cite{Liu,Tao,Wang2,Wang1},
where the nonlinear diffusions  in these studies were
still assumed to be non-degenerate.
For tumour angiogenesis model and relevant mathematical
analysis with or without degenerate diffusion, we refer to
\cite{Whang-zhian4,Mimura,Whang-zhian1,Whang-zhian2,Whang-zhian3,Winkler5} and the references therein.

Biological experiments suggest that
no cell migration (in particular no diffusivity) occurs in regions
where the tissue is absent \cite{Zhigun}.
In order to account for this biological feature,
we extend Chaplain and Anderson's
model \cite{Chaplain} to a new one with density-dependent jump probability
and quorum-sensing mechanisms of tumour cells as follows,
which is concerned with the competition between the following
several biological mechanisms: degenerate diffusion,
density-dependent chemotaxis, and logistic growth. That is,
\begin{equation}\label{eq-model-intro}
\begin{cases}
\displaystyle
\pd{u}{t}=\Delta (q(u)u)-\nabla\cdot(\phi(u)q(u)u\nabla v)+\mu u^{\delta}(1-ru),
\qquad &x\in \Omega,~t>0, \\[3mm]
\displaystyle
\pd{v}{t}=\Delta v+wz,
\qquad &x\in \Omega,~t>0,\\[3mm]
\displaystyle
\pd wt=-wz,
\qquad &x\in \Omega,~t>0,\\[3mm]
\displaystyle
\pd zt=\Delta z-z+u,
\qquad &x\in \Omega,~t>0.
\end{cases}
\end{equation}
The detailed derivation of the model \eqref{eq-model-intro} will be carried out in the last part of the next section.
Here,  $\Omega$ is a bounded domain in $\mathbb{R}^n$ with smooth boundary.
 The four variables $u$, $w$, $z$ and $v$ represent
the cancer cell density,
ECM concentration, the MDE concentration and the ECM* concentration, respectively.
 $q(u)$ denotes the jump probability of a cell depending on the population
pressure at its present location, which is increasing with respect to $u$ with
$q(0)=0$, $q(1)=1$, namely,
the jump probability is $1$ when the cell density exceeds maximum
and it is zero when the cell density is zero,
and $f(u)=\mu u^{\delta}(1-ru)$ is the logistic growth term,
where $\mu>0$ and $r>0$ are
the proliferation rate
and reciprocal of carrying capacity, respectively, $\delta\ge 1$ is a constant.
$\phi(u)$ is the density-dependent chemotactic
functions responding to quorum-sensing mechanisms, satisfying $|\phi(s)|\le 1$ and $|\phi'(s)|\le 1$.
While $\phi(u)$ can be sign-changing representing
the phenomenon that some chemicals
have been shown to elicit both attractive
and repellent responses \cite{Ming,Song}.
Moreover, some reasonable structure conditions on $\phi(s)$,
and $q(s)$ are also required in discussing the existence of solutions,
which we leave in Section 2 after the formulation of this model.

A recent interesting work related to the chemotaxis model mentioned
above is \cite{Ito},
in which they considered the following chemotaxis system with linear diffusion
\begin{equation}\label{eq-model-chemotaxis}
\begin{cases}
\displaystyle
\pd ut=\Delta u-\nabla\cdot(u\nabla v),
\qquad &x\in \Omega,~t>0, \\[2mm]
\displaystyle
\pd vt=\Delta v+wz,
\qquad &x\in \Omega,~t>0,\\[2mm]
\displaystyle
\pd wt=-wz,
\qquad &x\in \Omega,~t>0,\\[2mm]
\displaystyle
\pd zt=\Delta z-z+u,
\qquad &x\in \Omega,~t>0.
\end{cases}
\end{equation}
It is proved the existence of global solutions and the asymptotic
behaviors of global solutions as time goes to infinity
by using the properties of the Neumann heat semigroup $e^{t\Delta}$
in $\Omega$.

Compared to the linear cases,
the chemotactic system with degenerate diffusion
and chemotactic sensitivity is more complex and challenging.
Since the first equation of \eqref{eq-model-intro} is degenerate
at any point where $u(x,t)=0$,
there is no classical solution in general.
The spatial derivatives of $u$ may not exist in classical sense,
and may even do not belong to the class of locally integrable generalized functions,
that is, there might hold $u\not\in W^{2,1}_\text{loc}$.

As we all know,
for random walk equations, $u(x,t)>0$ for $t>0$ and any $x\in R^N$,
thus a linear diffusion process predicts a non-zero
number of the tumour cells for arbitrarily large displacements at arbitrarily small time,
namely, the underlying propagation speed is infinite \cite{Othmer}.
This means that
these models are valid for large time and describe the dynamics of tumour cells when they  spread to all parts of body.
Since most cancers develop at one anatomical site as primary tumour and then go on to metastasis,
it is vital to study the mechanism of tumour cell migration at the early stage
of cancer development when tumour progression begins and proceed
to yield a cancerous mass \cite{Weinberg}.

In this paper,
we provide a more realistic description of cell
migration process for early and late stages.
It is worth to mention that
our stability results of the model \eqref{eq-model-intro}
give a certain estimate for the expanding speed of tumour region before
cancer cells spread to the whole body.
We prove that there exist $t_0$ and
two families of monotone increasing open sets
$\{A_1(t)\}_{t>0}$, $\{A_2(t)\}_{t\in(0,t_0)}$ such that
$$A_1(t)\subset\text{supp }\ u(\cdot,t)\subset\overline A_2(t)\subset\Omega,
\quad t\in(0,t_0),$$
$\partial A_1(t)$ and $\partial A_2(t)$
have finite derivatives with respect to $t$,
namely, $\{A_1(t)\}_{t>0}$ and $\{A_2(t)\}_{t\in(0,t_0)}$ both expand
at finite speeds.
This indicates the finite speed propagation property of our chemotaxis model.
As shown late in Remarks \ref{remark2.1} and \ref{remark2.2}, in the porous media diffusion case,
we estimate that, in the early stage the expanding speed of tumour region is somehow like the algebraic rate of $(1+t)^\beta$ for some
$\beta \in (0,\frac{1}{2})$.

In contrast with the well known linear cases,
the degenerate diffusion is endowed with an interesting feature of slow diffusion,
that is, the compact support of solutions propagates at a finite speed.
The slow diffusion feature has some advantages and accuracy for describing
specified biological processes in the point of view of the physical reality,
and it also leads to more challenges in the mathematical studies.
For example, in order to investigate the asymptotic behavior of solutions,
one must appropriately describe the propagation behavior of its support,
which is more likely to be a compact subset of the prescribed domain
for some time interval if the initial data is given so.
We mention that the Neumann heat semigroup theory and functional transform methods
have been proved to be effective in studying the global boundedness and
large time behavior for the linear diffusion equations,
but they are all inapplicable in the degenerate diffusion case due to the nonlinearity.
We establish the global existence of bounded weak solutions to this model
by energy estimate technique and methods based on Moser-type iteration.
Then we prove that, as the late stage of the tumour migration, the original weak solution time-asymptotically
converges to its steady state, even if the initial perturbation is large, namely, the global stability of the steady state.
The adopted approach is the technical compactness analysis with the help of the comparison principle deduced by the approximate Hohmgren's approach
and two kinds of lower solutions
showing the expanding support and the exponentially convergence.
The one is a self similar weak lower solution of Barenblatt type
and the other kind is an ODE solution.

This paper is organized as follows.
In Section 2, we derive the models based on realistic biological assumptions,
which incorporate density-dependent jump probability
and quorum-sensing mechanisms and leads to the form of
equation given by \eqref{eq-model-intro}.
We leave the global existence of weak solutions and their regularity
to the corresponding chemotaxis system into Section 3 as preliminaries.
Section 4 is devoted to the study of compact support property of the tumour cells
at early stage and the large time behavior at late stage,
showing the exponentially convergence of solutions.

\section{Main results and formulations of new chemotaxis model}
In this section,
we first state our main results on the study of expanding compact support
of the tumour cells
at early stage and the asymptotic behavior at late stage. We leave the detailed derivation on the new chemotaxis model
\eqref{eq-model-intro} with
density-dependent jump probability and quorum-sensing mechanisms in the second part of this section.

\subsection{Main results}

We estimate the upper bound and lower bound for expanding speed of tumour cell region at early stage (before the tumour cells spread to the whole body) and
show the exponentially convergence of solutions for large time.
The derivation of this new chemotaxis model with degenerate diffusion
is presented at the end of this section.

We consider the following system \eqref{eq-model} with degenerate diffusion
\begin{equation}\label{eq-model}
\begin{cases}
u_t=\Delta (q(u)u)-\nabla\cdot(\phi(u)q(u)u\nabla v)+\mu u^{\delta}(1-u),\\
v_t=\Delta v+wz,\\
w_t=-wz,\\
z_t=\Delta z-z+u,
\quad &x\in \Omega,~t>0,\\[1mm]
\displaystyle
\pd{u}{\nu}=\pd{v}{\nu}=\pd{z}{\nu}=0,
\quad &x\in \partial\Omega,~t>0,\\
u(x,0)=u_0(x),\quad v(x,0)=v_0(x),\\
w(x,0)=w_0(x),\quad z(x,0)=z_0(x),
\quad &x\in\Omega,
\end{cases}
\end{equation}
where $\delta\ge1$,
$\mu>0$, $u_0,v_0,w_0,z_0$ are nonnegative functions,
$\nu$ is the unit outer normal vector,
and $q(u)\ge 0$ with $q(0)=0$. Here and after, the IBVP \eqref{eq-model} will be our main target equations.

Since degenerate diffusion equations may not have
classical solutions in general, we need to formulate the
following definition of generalized solutions
for the initial boundary value problem \eqref{eq-model}.

\begin{definition}
Let $T\in(0,\infty)$.
A quadruple $(u,v,w,z)$ is said to be a weak solution
to the problem \eqref{eq-model} in $Q_T=\Omega\times(0,T)$ if

{\rm(1)} $u\in L^\infty(Q_T)$, $\nabla (q(u)u)\in L^2((0,T);L^2(\Omega))$,
and $q(u)u_t\in L^2((0,T);L^2(\Omega))$;

{\rm(2)} $v\in L^\infty(Q_T)\cap L^2((0,T);W^{2,2}(\Omega))%
\cap W^{1,2}((0,T);L^2(\Omega))$;

{\rm(3)} $w\in L^\infty(Q_T)$, $w_t\in L^2((0,T);L^2(\Omega))$;

{\rm(4)} $z\in L^\infty(Q_T)\cap L^2((0,T);W^{1,2}(\Omega))%
\cap W^{1,2}((0,T);L^2(\Omega))$;

{\rm(5)} the identities
\begin{align*}
\int_0^T\int_\Omega u\psi_tdxdt&
+\int_\Omega u_0(x)\psi(x,0) dx
=\int_0^T\int_\Omega \nabla (q(u)u)\cdot\nabla\psi dxdt\\
&\qquad-\int_0^T\int_\Omega \phi(u)q(u)u\nabla v\cdot\nabla\psi dxdt
-\int_0^T\int_\Omega\mu u^\delta(1-u)\psi dxdt, \\
\int_0^T\int_\Omega v_t\varphi dxdt&
+\int_0^T\int_\Omega \nabla v\cdot\nabla\varphi dxdt
=\int_0^T\int_\Omega wz\varphi dxdt, \\
\int_0^T\int_\Omega w_t\varphi dxdt&
=-\int_0^T\int_\Omega wz\varphi dxdt, \\
\int_0^T\int_\Omega z_t\varphi dxdt&
+\int_0^T\int_\Omega \nabla z\cdot\nabla\varphi dxdt
=\int_0^T\int_\Omega (u-z)\varphi dxdt,
\end{align*}
hold for all $\psi,\varphi\in L^2((0,T);W^{1,2}(\Omega))%
\cap W^{1,2}((0,T);L^2(\Omega))$ with $\psi(x,T)=0$ for $x\in\Omega$;

{\rm(6)} $(v,w,z)$ takes the value $(v_0,w_0,z_0)$ in the sense of trace at $t=0$.

If $(u,v,w,z)$ is a weak solution of \eqref{eq-model} in $Q_T$ for
any $T\in(0,\infty)$, then we call it a global weak solution.

A quadruple $(u,v,w,z)$ is said to be a globally bounded weak solution
to the problem \eqref{eq-model} if there exists a constant $C$ such that
\begin{equation*}
\sup_{t\in\mathbb R^+}\left\{\|u\|_{L^\infty(\Omega)}+\|v\|_{W^{1,\infty}(\Omega)}
+\|w\|_{L^{\infty}(\Omega)}+\|z\|_{W^{1,\infty}(\Omega)}\right\}
\le C.
\end{equation*}
\end{definition}

Throughout this paper we assume that
$q(u)=u^{m-1}$ with $m>1$,
$|\phi(s)|\le1$, $|\phi'(s)|\le1$,
and the initial data satisfy
$u_0\in C^0(\overline\Omega)$,
$v_0\in W^{2,\infty}(\Omega)$,
$w_0\in C^{2,\theta}(\overline\Omega)$, $\theta\in(0,1)$,
$\pd{w_0}{\nu}=0$ on $\partial\Omega$,
$z_0\in C^0(\overline\Omega)$.
Here we note that for constant initial data $(u_0,v_0,w_0,z_0)$,
the first equation of \eqref{eq-model} is reduced to
$$u'(t)=\mu u^\delta(1-u), \quad u(0)=u_0,$$
which is ill-posed if $0<\delta<1$.
Therefore, we only consider the case $\delta\ge1$.

As preliminaries, we leave the global existence and regularity results into Section 3.
Our main results concerned with the description of cell migration process
are as follows.
First, we show that the tumour cells exist only in finite part of the body
at the early stage.

\begin{theorem}[Early stage profile - upper bound] \label{th-upper}
Let $(u,v,w,z)$ be a globally bounded weak solution of \eqref{eq-model}
with the initial data
$$\text{supp }u_0\subset\overline B_{r_0}(x_0)\subset\Omega,$$
for some $r_0>0$ and $x_0\in\Omega$.
Then there exists a time $t_1>0$ and a family of monotone increasing open sets
$\{A(t)\}_{t\in(0,t_1)}$ such that
$$\text{supp }u(\cdot,t)\subset\overline A(t)\subset\Omega, \quad t\in(0,t_1),$$
and $\partial A(t)$ has a finite derivative with respect to $t$.
More precisely, we can choose
$$A(t)=\{x\in\Omega; |x-x_0|^2<\eta(\tau+t)\}, \quad t\in(0,t_1),$$
with some appropriate $\eta, \tau>0$.
\end{theorem}

\begin{remark}\label{remark2.1}
As a typical finite propagating model,
the Barenblatt solution of the porous medium equation is
\begin{equation} \label{eq-Barenblatt}
B(x,t)=(1+t)^{-k}\Big[\Big(1-\frac{k(m-1)}{2mn}\frac{|x|^2}{(1+t)^{2k/n}}
\Big)_+\Big]^\frac{1}{m-1}
\end{equation}
with $k=1/(m-1+2/n)<n/2$ for $m>1$,
and its support is expanding at the rate $(1+t)^{k/n}$.
Here we have proved the tumour cells are located within a ball
expanding at the rate $(1+t)^{1/2}$.
\end{remark}

Next, we show the propagating property of the tumour cells at the early stage.

\begin{theorem}[Early stage profile - lower bound] \label{th-lower}
Let $(u,v,w,z)$ be a globally bounded weak solution of \eqref{eq-model}.
Assume that $1\le\delta<m$, $\Omega$ is convex and $u_0\not\equiv0$.
Then there exists a time $t_2>0$ such that the support of $u$ expands
to the whole $\Omega$ when $t\ge t_2$.
Precisely speaking, there exist a family of monotone increasing open sets
$\{A(t)\}_{t>0}$
(we can choose $A(t)=\{x\in\Omega; |x-x_0|^2<\eta(1+t)^\beta\}$
with sufficiently small $\beta, \eta>0$)
such that
$$A(t)\subset\text{supp }u(\cdot,t), \quad t>0,$$
and $A(t)=\Omega$ for $t\ge t_2$,
$\partial A(t)$ has a finite derivative with respect to $t$.
\end{theorem}

\begin{remark}\label{remark2.2}
For this chemotaxis system, we proved that
the tumour cells will expand to the whole body when the time $t$ increases.
Compared with the porous medium equation,
whose Barenblatt solution $B(x,t)$ in \eqref{eq-Barenblatt}
is expanding at the rate $(1+t)^{2k/n}$,
the tumour cells of \eqref{eq-model} migrate to at least a ball expanding at the rate $(1+t)^\beta$.
Here in the proof we have selected $\beta>0$ sufficiently small,
which means the support is expanding with a much slower rate.
\end{remark}

Under the hypotheses of Theorem \ref{th-upper} and Theorem \ref{th-lower},
we see that there exist $t_0$ and
two family of monotone increasing open sets
$\{A_1(t)\}_{t>0}$, $\{A_2(t)\}_{t\in(0,t_0)}$ such that
$$A_1(t)\subset\text{supp }u(\cdot,t)\subset\overline A_2(t)\subset\Omega,
\quad t\in(0,t_0),$$
$\partial A_1(t)$ and $\partial A_2(t)$
have finite derivatives with respect to $t$,
which means that $\{A_1(t)\}_{t>0}$ and $\{A_2(t)\}_{t\in(0,t_0)}$ both expand
at finite speeds.
This indicates immediately the finite speed propagation property of this chemotaxis model,
though we have not proved it directly.

After the tumour cells spread to the whole body, we can investigate the large time behavior.
We show that the solution converges to its steady state exponentially.

\begin{theorem}[Late stage profile] \label{th-asymp}
Let $(u,v,w,z)$ be a globally bounded weak solution of \eqref{eq-model}.
Assume that the hypothesis in Theorem \ref{th-lower} is valid.
Then there exist $C$ and $c>0$ such that
$$\|u(\cdot,t)-1\|_{L^{\infty}(\Omega)}
+\|w(\cdot,t)\|_{W^{1,\infty}(\Omega)}
+\|v(\cdot,t)-(\overline v_0+\overline w_0)\|_{W^{2,\infty}(\Omega)}
+\|z(\cdot,t)-1\|_{L^{\infty}(\Omega)}
\le Ce^{-ct},$$
for all $t>0$, where
$\overline v_0=\frac{1}{|\Omega|}\int_\Omega v_0(x) dx$
and $\overline w_0=\frac{1}{|\Omega|}\int_\Omega w_0(x) dx$.
\end{theorem}

The main difficulty lies in proving the expanding property of the support
of the first component.
We first prove the comparison principle by the approximate Hohmgren's approach,
and then construct two kinds of lower solutions.
The one is a self similar weak lower solution with
much slower expanding support and slightly faster decaying maximum
compared with the Barenblatt solution to the porous medium equation,
 the other kind is an ODE solution.
After showing the expanding property, we formulate several upper and lower solutions
that converge to steady state exponentially by utilizing the exponential decay of
other components.

\subsection{Derivation of the new chemotaxis model}
We extend the derivation of the classical
 taxis models in \cite{Othmer}.
The derivation of the model begins with a master equation
for a continuous-time and discrete-space random walk
\begin{equation}\label{eq-Master}
\pd{u_i}{t}=\mathcal T_{i-1}^+u_{i-1}
+\mathcal T_{i+1}^-u_{i+1}-(\mathcal T_i^++\mathcal T_i^-)u_i,
\end{equation}
where $\mathcal T^{\pm}_i(\cdot)$ denote the transitional-probabilities
per unit time of a one-step jump to $i\pm1$
and $u_i$ denotes the cell density at $i$.

Cancer cells can modify their
migration mechanisms in response to different conditions\cite{Friedl}.
There are two potentially important factors:
(i) the effect of cell-density on the probability of cell movement;
(ii) the effect of signal-mediated cell-density sensing mechanisms
on movement \cite{Painter}.

For neighbor-based and gradient-based rules,
Painter and Hillen \cite{Painter} proposed volume filling approach,
that is, the movement of cells is inhibited by the neighboring site
where the cells are densely packed.
The transitional probability then takes the form
\begin{equation}\label{eq-T_Painter}
\mathcal T^{\pm}_i=q(u_{i\pm1})(\alpha
+\beta(\tau(v_{i\pm1})-\tau(v_i))),
\end{equation}
where $q(u)$ denotes the probability of a cell finding
space at its neighboring location,
constant $\alpha$ is the intrinsic dispersion coefficient,
constant $\beta$ the coefficient signal detection,
$v_i$ the signal concentration,
and $\tau$ the mechanism of tactic responses
in cell populations, such as chemotaxis, haptotaxis or phototaxis.
Substituting \eqref{eq-T_Painter} to the master equation \eqref{eq-Master},
in the PDE limits they derives
\begin{equation*}
\pd u t=\nabla\cdot(d_1(q(u)-q'(u))\nabla u-\chi(v)q(u)u\nabla u)
\end{equation*}
where $d_1=k\alpha$, $\chi(v)=2k\beta\frac{d\tau(v)}{dv}$,
$k$ is a scaling constant.
Note that $q(u)$ is a non-increasing function in this model,
which says that the probability of a cell finding space at its
neighboring site decreases in the cell density at that site.

Since a different combination of the above strategies
may be necessary to reflect cell movement,
we combine the local and gradient-based strategies and assume
the transitional probability of the form
\begin{equation}\label{eq-T}
\mathcal T^{\pm}_i=q(u_i)(\alpha
+\beta(\tau(v_{i\pm1})-\tau(v_i))),
\end{equation}
where $q(u)$ represents the jump probability of a cell
due to the population pressure at present site.
At the microscopic level,
a high cell density results in increased probability of a cell being
``pushed'' from departure site \cite{Lou,Okubo,Sherratt},
for example due to the pressure exerted by neighboring cells.
We shall assume that only a finite number of cells,
$U_\text{max}$, can be accommodated at any site.
We study the relative density $\tilde u=u/U_\text{max}$,
(and drop the symbol $\tilde{}$ for simplicity).
Moreover, the jump probability is $1$ when the cell
density exceeds $U_\text{max}$
and it is zero when the cell density is zero.
Thus we stipulate the following conditions on $q$:
$$
q(0)=0, \quad q(1)=1 \quad \text{and}\quad q(u)\ge0, \quad
\text{for all~}0\le u\le1.
$$
A natural choice for $q(u)$ is
\begin{equation} \label{eq-qu}
q(u)=u^{m-1},\quad m>1,
\end{equation}
which states that the probability of a jump leaving one site increases
with the cell density at that site \cite{Murry,Rodrigo}.

In most tissues, to control cell density at proper level,
cells also secrete quorum-sensing molecule $z$ ,
then the concentration of the modules sensed by the cell will
be an indicator of local density \cite{Painter}.
Cells sense the same gradient of chemical at the surface,
but the ``strength'' signalled to the movement dynamics is
modulated by quorum-sensing molecule $z$ \cite{Song}.
A more general choice of transitional-probabilities
$\mathcal T^{\pm}_i(\cdot)$ can also be considered, namely
\begin{equation}\label{Tbetau}
\mathcal T^{\pm}_i=q(u_i)(\alpha
+\beta(z_i)(\tau(v_{i\pm1})-\tau(v_i))),
\end{equation}
where $\beta(z)$ is a tactic function
responding to quorum-sensing molecule $z$.
Substituting \eqref{Tbetau} into the Master Equation \eqref{eq-Master} gives:
\begin{align*}
\frac{d}{dt}u_i=&
q_{i-1}(\alpha+\beta_{i-1}(\tau_i-\tau_{i-1}))u_{i-1}
+
q_{i+1}(\alpha+\beta_{i+1}(\tau_i-\tau_{i+1}))u_{i+1}
\\&-
q_{i}(\alpha+\beta_i(\tau_{i+1}-\tau_{i}))u_{i}
-
q_{i}(\alpha+\beta_i(\tau_{i-1}-\tau_{i}))u_{i}
\\
=&\alpha(q_{i-1}u_{i-1}+q_{i+1}u_{i+1}-2q_iu_i)
\\&
+\beta_{i-1}
q_{i-1}(\tau_i-\tau_{i-1})u_{i-1}
+\beta_{i+1}q_{i+1}(\tau_i-\tau_{i+1})u_{i+1}
-\beta_iq_i(\tau_{i+1}+\tau_{i-1}-2\tau_i)u_i\\
=&\alpha(q_{i-1}u_{i-1}+q_{i+1}u_{i+1}-2q_iu_i)
\\
&-\beta_{i+1}
q_{i+1}u_{i+1}(\tau_{i+1}-\tau_i)
+\beta_iq_iu_i(\tau_i-\tau_{i-1})
-
\big(\beta_i
q_iu_i(\tau_{i+1}-\tau_i)
-\beta_{i-1}q_{i-1}u_{i-1}(\tau_{i}-\tau_{i-1})
\big)\\
=&\alpha(q_{i-1}u_{i-1}+q_{i+1}u_{i+1}-2q_iu_i)
\\
&-\Big(
(\beta_{i+1}q_{i+1}u_{i+1}+\beta_iq_iu_i)(\tau_{i+1}-\tau_i)
-(\beta_{i-1}q_{i-1}u_{i-1}+\beta_iq_iu_i)(\tau_{i}-\tau_{i-1})
\Big).
\end{align*}
We set $x = kh$, interpret $x$ as a continuous variable and extend the
definition of $u_i$ accordingly.
The transitional probabilities of jumping to a neighboring location
depend on the spatial scale $h$.
Thus we assume that
$\mathcal T^{\pm}_h=\frac{k}{h^2}\mathcal T^{\pm}$
for some scaling constant $k$.
Expanding the right-hand side with respect to $h$,
we obtain for the cell density $u(x,t)$:
$$
\pd{u}{t}=k\Big(
\alpha\frac{\partial^2(q(u)u)}{\partial x^2}
-2\frac{\partial}{\partial x}\Big(\beta(z)q(u)u\pd{\tau}{x}
\Big)
\Big)+O(h^2).
$$
By taking the limit of $h\to0$, we arrive at the following model
$$
\pd{u}{t}=
D_u\frac{\partial^2(q(u)u)}{\partial x^2}
-\frac{\partial}{\partial x}\Big(\beta(z)\chi(v)q(u)u\pd{v}{x}
\Big),
$$
where $D_u=k\alpha$, $\chi(v)=2k\frac{d\tau(v)}{dv}$.
The function $\chi(v)$ is commonly referred as
the tactic sensitivity function.
The simplest form  is $\chi(v)=\chi_0$
with $\chi_0$ being a constant.

Apart from that,
we consider a modification of the Verhulst
logistic growth term to model organ size evolution
introduced by Blumberg \cite{Blumberg} and Turner \cite{Turner},
which is called hyper-logistic function, accordingly
$$
f(u)=r u^\delta(1-\mu u).
$$
In the special case,  the quorum sensing molecule $z=z(u)$ is
not diffusing and a monotone increasing function of the
cell density.
Denote $\beta(z)=\beta(z(u)):=\phi(u)$.
Assume that $z$ switches the response to chemotaxis concentration $v$
from attractant at low concentrations of $v$
to repellent at high concentrations,
namely, $\beta$ is a sign-changing and non-increasing function.
(e.g. $\beta(z)=1-{z}/{z*}$) \cite{Painter, Painter19},
Including cell kinetics and signal dynamics,
we derive the resulting model for the cell movement
$$
\pd ut=\underbrace{D_u\Delta(q(u)u)}_{\text{dispersion}}
-\underbrace{\chi_0
\nabla\cdot(\phi(u)q(u)u\nabla v)}_{\text{chemotaxis}}
+\underbrace{\mu u^\delta(1-r u)}_{\text{proliferation}}.
$$
Incorporating the kinetic equation of ECM and MDE,
we arrive at a modified Chaplain and Lolas' chemotaxis model,
see  \eqref{eq-model}, where we assume the constants $D_u, \chi_0, r=1$
for simplification.

\section{Preliminaries: Global existence, boundedness and regularity}

As preliminaries, we prove the existence, boundedness and regularity of a global weak solution in this section.
The main preliminary results are as follows.

\begin{theorem}[Existence of globally bounded weak solutions] \label{th-existence}
For $1\le n\le3$,
the problem \eqref{eq-model} admits a
globally bounded weak solution $(u,v,w,z)$.
\end{theorem}

\begin{theorem}[Regularity] \label{th-regular}
Let $(u,v,w,z)$ be a globally bounded weak solution of \eqref{eq-model}.
Then there exist $\alpha\in(0,1)$ and $C(p)>0$ such that
$$\|u\|_{L^\infty(\Omega\times(t,t+1))}
+\|v\|_{C^{2+\alpha,1+\alpha/2}(\overline\Omega\times[t,t+1])}
+\|w\|_{C^\alpha(\overline\Omega\times[t,t+1])}
+\|z\|_{W^{2,1}_{p}(\Omega\times(t,t+1))}\le C(p),$$
for any $p>1$ and $t\ge1$.
\end{theorem}

We first use the artificial viscosity method to get smooth approximate solutions.
Despite the absence of comparison principle,
we can prove a special case compared with a lower solution,
which is helpful for establishing the regularity estimates.
By making use of the special structure of dispersion, we carry on the estimates on
$u^m$ in $W^{1,2}(Q_T)$, instead of $u$.
These energy estimates ensure the global existence of weak solution.

Consider the following corresponding regularized problem
\begin{equation}\label{eq-model-regularized}
\begin{cases}
u_t=\nabla\cdot(m(a_\varepsilon(u))^{m-1}\nabla u)-\nabla\cdot(u^m\phi(u)\nabla v)
+\mu|u|^{\delta-1}u(1-u)+\varepsilon,\\
v_t=\Delta v+wz,\\
w_t=-wz,\\
z_t=\Delta z-z+u,
\quad &x\in \Omega,~t>0,\\[2mm]
\displaystyle
\pd{u}{\nu}=\pd{v}{\nu}=\pd{z}{\nu}=0,
\quad &x\in \partial\Omega,~t>0,\\
u(x,0)=u_{0\varepsilon}(x),\quad v(x,0)=v_{0\varepsilon}(x),\\
w(x,0)=w_{0\varepsilon}(x),\quad z(x,0)=z_{0\varepsilon}(x),
\quad &x\in\Omega,
\end{cases}
\end{equation}
where $\varepsilon\in(0,1)$, $a_\varepsilon\in C^\infty(\mathbb R)$,
$a_\varepsilon(s)=s+\varepsilon$ for $s\ge0$,
$a_\varepsilon(s)=\varepsilon/2$ for $s<-\varepsilon$, $a_\varepsilon$ is
monotone increasing with $0\le a'_\varepsilon\le1$,
and $u_{0\varepsilon},v_{0\varepsilon},w_{0\varepsilon},z_{0\varepsilon}$
are smooth approximations of $u_0,v_0,w_0,z_0$, respectively, with
\begin{align*}
\varepsilon\le u_{0\varepsilon}\le u_0+\varepsilon,& \quad
0\le v_{0\varepsilon}\le v_0+\varepsilon, \\
0\le w_{0\varepsilon}\le w_0+\varepsilon,& \quad
0\le z_{0\varepsilon}\le z_0+\varepsilon, \\
|\nabla u_{0\varepsilon}|\le2|\nabla u_0|,& \quad
|\nabla v_{0\varepsilon}|\le2|\nabla v_0|, \\
|\nabla w_{0\varepsilon}|\le2|\nabla w_0|, \quad
|\Delta w_{0\varepsilon}|&\le2|\Delta w_0|, \quad
|\nabla z_{0\varepsilon}|\le2|\nabla z_0|,
\end{align*}
and $\pd{w_{0\varepsilon}}{\nu}=0$ on $\partial\Omega$.
The local existence of the regularized problem \eqref{eq-model-regularized}
is trivial and we denote the unique solution
by $(u_\varepsilon,v_\varepsilon,w_\varepsilon,z_\varepsilon)$.
Let $(0,T_{\text{max}})$ be its maximal existence interval.

As usual, there is no comparison principle for the system, because the system is strongly coupled.
However, we have the following lemma.

\begin{lemma} \label{le-comparison}
There holds $u_\varepsilon\ge0$, $v_\varepsilon\ge0$, $w_\varepsilon\ge0$, and
$z_\varepsilon\ge0$ for all $x\in\Omega$ and $t\in(0,T_{\text{max}})$.
\end{lemma}
{\it \bfseries Proof.}
We denote $(u_\varepsilon,v_\varepsilon,w_\varepsilon,z_\varepsilon)$
by $(u,v,w,z)$ in this proof for the sake of simplicity.
We argue by contradictions.
Since $u_{0\varepsilon}\ge\varepsilon>0$, there exists $t_0\in(0,T_{\text{max}})$ such that
$u>0$ for all $x\in\Omega$ and $t\in(0,t_0)$,
$u(x_0,t_0)=0$ for some $x_0\in\overline\Omega$ and $u(x,t_0)\ge0$ for all $x\in\Omega$.

Now we divide this proof into two parts.
If $x_0\in\Omega$, then $\nabla u(x_0,t_0)=0$ and
\begin{align*}
\nabla\cdot(m(a_\varepsilon(u))^{m-1}\nabla u)
&=m(a_\varepsilon(u))^{m-1}\Delta u+m(m-1)a'_\varepsilon(u)|\nabla u|^2\ge0, \\
\nabla\cdot(u^m\phi(u)\nabla v)&=u^m\phi(u)\Delta v
+(mu^{m-1}\phi(u)+u^m\phi'(u))\nabla u\cdot\nabla v=0, \\
\mu |u|^{\delta-1}u&(1-u-w)=0,
\end{align*}
which contradicts to $\pd{u}{t}(x_0,t_0)\le0$.

If $x_0\in\partial\Omega$, then $\pd{u}{\tau}(x_0,t_0)=0$,
$\frac{\partial^2u}{\partial \tau^2}(x_0,t_0)\ge0$
for any tangent vector $\tau$, and
the boundary condition shows that $\pd{u}{\nu}(x_0,t_0)=0$.
We assert that $\frac{\partial^2u}{\partial \nu^2}(x_0,t_0)\ge0$.
In fact, if it were not true, Taylor expansion at $(x_0,t_0)$ shows that
there would exist a point $x'\in\Omega$ such that $u(x',t_0)<0$.
Therefore, we also have $\nabla u(x_0,t_0)=0$ and the above equalities.
Those contradictions complete the proof.
$\hfill\Box$

Since $u_\varepsilon\ge0$, the first equation of \eqref{eq-model-regularized}
is equivalent to
\begin{equation*}
\pd{u}{t}=\Delta(u+\varepsilon)^m-\nabla\cdot(u^m\phi(u)\nabla v)
+\mu u^{\delta}(1-u)+\varepsilon, \qquad u\ge0.
\end{equation*}

Now we present some energy estimates
independent of time $t$ and the parameter $\varepsilon$.

\begin{lemma} \label{le-uL1}
The first solution component $u_\varepsilon$ satisfies
\begin{equation*}
\sup_{t\in(0,T_{\text{max}})}\int_\Omega u_\varepsilon(\cdot, t)dx\le
\max \left\{\int_\Omega u_0 dx+|\Omega|,
\left(\frac{2(C_1+|\Omega|)}{\mu C_2}\right)^{1/(\delta+1)}\right\},
\end{equation*}
where $C_1=\mu2^\delta|\Omega|$ and $C_2=1/|\Omega|^{\delta}$.
\end{lemma}
{\it \bfseries Proof.}
We denote $u_\varepsilon$ by $u$ in this proof for the sake of simplicity.
Since $u$ is nonnegative and $\pd{u}{\nu}=\pd{v}{\nu}=0$ on $\partial \Omega$,
integration of the first equation of \eqref{eq-model-regularized} over $\Omega$ yields
\begin{equation}\nonumber
\frac{d}{dt}\int_\Omega u dx\le\mu\int_\Omega u^\delta dx
-\mu\int_\Omega u^{\delta+1}dx+|\Omega|,
\end{equation}
for all $t\in (0,T_{\text{max}})$.
We note that
$$\mu\int_\Omega u^\delta dx\le\frac{1}{2}\mu\int_\Omega u^{\delta+1}dx+C_1,$$
and
$$\int_\Omega u^{\delta+1}dx\ge C_2\left(\int_\Omega u dx\right)^{\delta+1},$$
where $C_1=\mu2^\delta|\Omega|$ and $C_2=1/|\Omega|^{\delta}$.
Let $y(t)=\int_\Omega u(\cdot, t)dx$ for $t\in[0,T_{\text{max}})$.
We find
\begin{equation}\nonumber
y'(t)\le C_1+|\Omega|-\frac{\mu C_2}{2} y^{\delta+1}(t).
\end{equation}
The  comparison principle of ODE shows that
\begin{equation}\nonumber
y(t)\le\max\left\{y(0),\left(\frac{2(C_1+|\Omega|)}{\mu C_2}\right)^{1/{\delta+1}}\right\}
\end{equation}
for all $t\in (0,T_{\text{max}})$.
$\hfill\Box$

Here we recall some lemmas about the $L^p$-$L^q$ type estimates
for the components of the solution,
and we refer the readers to \cite{Ito} for details.

\begin{lemma}[\cite{Ito}] \label{le-zLq}
Let $p\ge1$ and
\begin{align*} \label{eq-zLq}
\begin{cases}
q\in[1,\frac{np}{n-2p}), \qquad &p\le\frac{n}{2},\\
q\in[1,\infty], \qquad &p>\frac{n}{2}.
\end{cases}
\end{align*}
Then for any $T\in(0,T_\text{max}]$, there exists a constant $C_z(p,q)$ such that
\begin{equation*}
\sup_{t\in(0, T)}\|z_\varepsilon(\cdot,t)\|_{L^q(\Omega)}
\le C_z(p,q)(\|z_0\|_{L^q(\Omega)}
+\sup_{t\in(0, T)}\|u_\varepsilon(\cdot,t)\|_{L^p(\Omega)}).
\end{equation*}
\end{lemma}

\begin{lemma}[\cite{Ito}] \label{le-nablavLr}
Let $q\ge1$ and
\begin{align*}
\begin{cases}
r\in[1,\frac{nq}{n-q}), \qquad &q\le n,\\
r\in[1,\infty], \qquad &q>n.
\end{cases}
\end{align*}
Then for any $T\in(0,T_\text{max}]$, there exists a constant $C_v(q,r)$ such that
\begin{equation*}
\sup_{t\in(0, T)}\|\nabla v_\varepsilon (\cdot,t)\|_{L^r(\Omega)}
\le C_v(q,r)(\|\nabla v_0\|_{L^r(\Omega)}
+\sup_{t\in(0, T)}\|z_\varepsilon(\cdot,t)\|_{L^q(\Omega)}).
\end{equation*}
\end{lemma}

\begin{lemma} \label{le-wv}
There holds
$$\|w_\varepsilon(\cdot,t)\|_{L^\infty(\Omega)}
\le\|w_0\|_{L^\infty(\Omega)}+1,
\quad t\in(0,T_\text{max}),$$
and
$$\int_\Omega v_\varepsilon(x,t)dx
\le \int_\Omega v_0(x)dx+\int_\Omega w_0(x)dx+2|\Omega|,
\quad t\in(0,T_\text{max}).$$
\end{lemma}
{\it \bfseries Proof.}
Since both $w_\varepsilon$ and $z_\varepsilon$ are nonnegative,
it is clear from the third equation of \eqref{eq-model-regularized} that
$$|w_\varepsilon(x,t)|\le w_{0\varepsilon}(x,t)\le \|w_0\|_{L^\infty(\Omega)}+1.$$
We add the third to the second equation of \eqref{eq-model-regularized}
and integrate over $\Omega$ to obtain
$$\frac{d}{dt}\int_\Omega(v_\varepsilon+w_\varepsilon)dx
=\int_\Omega \Delta v_\varepsilon dx=0, \quad t\in(0,T_\text{max}).$$
Thus,
$$\int_\Omega (v_\varepsilon+w_\varepsilon)dx
\le \int_\Omega v_{0\varepsilon}(x)dx+\int_\Omega w_{0\varepsilon}(x)dx
\le \int_\Omega v_0(x)dx+\int_\Omega w_0(x)dx+2|\Omega|,$$
for all $t\in(0,T_\text{max})$.
$\hfill\Box$

\begin{lemma} \label{le-nablavLany}
Let $1\le n\le 3$.
There exists a constant $C$ independent of $t$ and $\varepsilon$ such that
$$\|v_\varepsilon\|_{L^\infty(\Omega)}\le C,
\quad t\in(0,T_\text{max}).$$
For any $r\ge1$, there exists a constant $C(r)$
independent of $t$ and $\varepsilon$ such that
$$\|\nabla v_\varepsilon\|_{L^r(\Omega)}\le C(r),
\quad t\in(0,T_\text{max}).$$
\end{lemma}
{\it \bfseries Proof.}
According to Lemma \ref{le-uL1}, $\|u_\varepsilon\|_{L^1(\Omega)}$ is uniformly bounded.
Since $n\le3$, we can apply Lemma \ref{le-zLq} and \ref{le-nablavLr}
to complete this proof.
$\hfill\Box$

The following Gagliardo-Nirenberg inequality (see \cite{Wang1,GN})
will be used in deriving the $L^p$ estimates of $u_\varepsilon$.

\begin{lemma} \label{le-GN}
Let $0<s\le p\le\frac{2n}{(n-2)_+}$.
There exists a positive constant $C$ such that for all
$u\in W^{1,2}(\Omega)\cap L^s(\Omega)$,
$$\|u\|_{L^p(\Omega)}\le C(\|\nabla u\|_{L^2(\Omega)}^a\|u\|_{L^s(\Omega)}^{1-a}
+\|u\|_{L^s(\Omega)})$$
is valid with $a=\frac{n/s-n/p}{1-n/2+n/s}\in(0,1)$.
\end{lemma}

We present the following $L^p$ estimate of $u_\varepsilon$.

\begin{lemma} \label{le-uLp}
Let $1\le n\le3$.
For any given $p\ge1$, there exists a constant $C(p)>0$
independent of $t$ and $\varepsilon$ such that
\begin{equation*}
\|u_\varepsilon(\cdot,t)\|_{L^p(\Omega)}\le C(p),
\quad t\in(0,T_\text{max}).
\end{equation*}
\end{lemma}
{\it \bfseries Proof.}
We denote $u_\varepsilon,v_\varepsilon$ by $u,v$ in this proof for the sake of simplicity.
By a straightforward computation, testing the
first equation in \eqref{eq-model-regularized} by $u^r$ for $r>0$
and integrating by parts we find that
\begin{align} \nonumber
&\frac{1}{r+1}\frac{d}{dt}\int_\Omega u^{r+1}dx
+\int_\Omega \nabla(u+\varepsilon)^m\cdot \nabla u^r dx
\\ \label{eq-zup}
\le&\int_\Omega u^m\phi(u)\nabla v\cdot \nabla u^r dx
+\mu\int_\Omega u^{\delta+r}dx-\mu\int_\Omega u^{\delta+r+1}dx
+\int_\Omega u^r dx.
\end{align}
We note that
\begin{equation} \label{eq-zuap1}
\mu\int_\Omega u^{\delta+r}dx\le\frac{1}{4}\mu\int_\Omega u^{\delta+r+1}dx+C_1,
\end{equation}
and
\begin{equation} \label{eq-zupm1}
\int_\Omega u^r dx\le\frac{1}{4}\mu\int_\Omega u^{\delta+r+1}dx+C_2,
\end{equation}
where $C_1$ and $C_2$ are constants independent of $t$.
Then by Young's inequality, we see that
\begin{align} \nonumber
\int_\Omega u^m\phi(u)\nabla v\cdot \nabla u^r dx
&\le r\int_\Omega u^{m+r-1}|\nabla v\cdot\nabla u|dx
\\ \nonumber
&\le\frac{mr}{2}\int_\Omega(u+\varepsilon)^{m-1}u^{r-1}|\nabla u|^2dx
+\frac{r}{2m}\int_\Omega u^{m+r}|\nabla v|^2dx
\\\label{eq-znablav}
&\le\frac{1}{2}\int_\Omega\nabla(u+\varepsilon)^m\cdot\nabla u^r dx
+\frac{r}{2m}\int_\Omega u^{m+r}|\nabla v|^2dx.
\end{align}
We use H\"older's inequality to see that
\begin{align*}
\frac{r}{2m}\int_\Omega u^{m+r}|\nabla v|^2dx
&\le\frac{r}{2m}\Big(\int_\Omega u^{m+r+\kappa}dx\Big)^{\frac{m+r}{m+r+\kappa}}
\Big(\int_\Omega |\nabla v|^{\frac{2(m+r+\kappa)}{\kappa}}dx\Big)%
^{\frac{\kappa}{m+r+\kappa}}\\
&\le C_3\Big(\int_\Omega u^{m+r+\kappa}dx\Big)^{\frac{m+r}{m+r+\kappa}}
\end{align*}
where $\kappa>0$ is a constant to be determined
and $C_3$ is a constant depending on the
$L^\frac{2(m+r+\kappa)}{\kappa}(\Omega)$ norm of $\nabla v$
which is uniformly bounded according to Lemma \ref{le-nablavLany}.
Now we use the Gagliardo-Nirenberg inequality Lemma \ref{le-GN} to obtain
\begin{align*}
\Big(\int_\Omega u^{m+r+\kappa}dx\Big)^{\frac{m+r}{m+r+\kappa}}
&=\|u^\frac{m+r}{2}\|_{L^{\frac{2(m+r+\kappa)}{m+r}}(\Omega)}^2\\
&\le C_4\big(\|\nabla u^\frac{m+r}{2}\|_{L^2(\Omega)}^{2a}%
\|u^\frac{m+r}{2}\|_{L^{\frac{2}{m+r}}(\Omega)}^{2(1-a)}
+\|u^\frac{m+r}{2}\|_{L^{\frac{2}{m+r}}(\Omega)}^2\big)\\
&\le C_5(1+\|\nabla u^\frac{m+r}{2}\|_{L^2(\Omega)}^{2a}),
\end{align*}
where $C_4$ is a constant, $C_5$ depends on $\|u\|_{L^1(\Omega)}$, and
$$a=\frac{n(m+r)/2-n(m+r)/(2(m+r+\kappa))}{1-n/2+n(m+r)/2}\in(0,1),$$
provided that
$\frac{2(m+r+\kappa)}{m+r}<\frac{2n}{(n-2)_+}$.
This can be done by taking $\kappa>0$ and appropriately small.
Therefore, we have
\begin{align}\nonumber
\frac{r}{2m}\int_\Omega u^{m+r}|\nabla v|^2dx
&\le C_3C_5(1+\|\nabla u^\frac{m+r}{2}\|_{L^2(\Omega)}^{2a})\\ \nonumber
&\le\frac{2mr}{(m+r)^2}\|\nabla u^\frac{m+r}{2}\|_{L^2(\Omega)}^2+C_6\\ \nonumber
&\le\frac{mr}{2}\int_\Omega(u+\varepsilon)^{m-1}u^{r-1}|\nabla u|^2dx+C_6
\\ \label{eq-zgn}
&\le\frac{1}{2}\int_\Omega\nabla(u+\varepsilon)^m\cdot\nabla u^r dx+C_6,
\end{align}
since $a\in(0,1)$.
Combining \eqref{eq-zuap1}, \eqref{eq-zupm1},
\eqref{eq-znablav}, \eqref{eq-zgn} with \eqref{eq-zup}, we infer that
\begin{equation}\nonumber
\frac{d}{dt}\int_\Omega u^{r+1}dx\le
-\frac{\mu(r+1)}{2}\int_\Omega u^{\delta+r+1}dx+(r+1)(C_1+C_2+C_6).
\end{equation}
According to
$$\int_\Omega u^{\delta+r+1}dx\ge
\frac{1}{|\Omega|^\frac{\delta}{r+1}}%
\Big(\int_\Omega u^{r+1}dx\Big)^{\frac{\delta+r+1}{r+1}},$$
we obtain
\begin{align*}
\frac{d}{dt}\int_\Omega u^{r+1}dx
\le (r+1)(C_1+C_2+C_6)-\frac{\mu(r+1)}{2|\Omega|^\frac{\delta}{r+1}}%
\Big(\int_\Omega u^{r+1}dx\Big)^{\frac{\delta+r+1}{r+1}}.
\end{align*}
By an ODE comparison,
\begin{equation*}
\int_\Omega u^{r+1}dx\le\max\left\{\int_\Omega (u_0+1)^{r+1}dx,
\Big(\frac{2(C_1+C_2+C_6)|\Omega|^\frac{\delta}{r+1}}{\mu}
\Big)^{\frac{r+1}{\delta+r+1}}\right\}
\end{equation*}
for all $t\in(0,T)$.
$\hfill\Box$

\begin{lemma} \label{le-nablavinfty}
Let $1\le n\le3$.
There exists a constant $C>0$ independent of $T_\mathrm{max}$
and $\varepsilon$ such that
\begin{equation*}
\sup_{t\in(0, T_\text{max})}\|\nabla v_\varepsilon\|_{L^\infty(\Omega)}\le C.
\end{equation*}
\end{lemma}
{\it \bfseries Proof.}
According to Lemma \ref{le-uLp}, $\|u_\varepsilon\|_{L^{n+1}(\Omega)}$
is uniformly bounded.
We can apply Lemma \ref{le-zLq} and Lemma \ref{le-nablavLr}
to obtain the boundedness of $\|\nabla v_\varepsilon\|_{L^\infty(\Omega)}$.
$\hfill\Box$

We now employ the following Moser-type iteration to get the
$L^\infty(\Omega)$ estimate of $u$.

\begin{lemma} \label{le-uLinfty}
Let $1\le n\le3$.
There exists a constant $C>0$ independent of $T_\mathrm{max}$
and $\varepsilon$ such that
\begin{equation*}
\sup_{t\in(0, T_\text{max})}\|u_\varepsilon\|_{L^\infty(\Omega)}\le C.
\end{equation*}
\end{lemma}
{\it \bfseries Proof.}
We denote $u_\varepsilon,v_\varepsilon$ by $u,v$ in this proof for the sake of simplicity.
We test the first equation in \eqref{eq-model-regularized} by $u^r$ for $r>0$
and integrating by parts we find that
\begin{align} \nonumber
&\frac{1}{r+1}\frac{d}{dt}\int_\Omega u^{r+1}dx
+\int_\Omega \nabla(u+\varepsilon)^m\cdot \nabla u^r dx
\\ \label{eq-zur}
\le&\int_\Omega u^m\phi(u)\nabla v\cdot \nabla u^r dx
+\mu\int_\Omega u^{\delta+r}dx-\mu\int_\Omega u^{\delta+r+1}dx
+\int_\Omega u^r dx.
\end{align}
Similar to the proof of Lemma \ref{le-uLp},
using Young's inequality we can estimate
\begin{align*}
\mu\int_\Omega u^{\delta+r}dx&\le\frac{1}{4}\mu\int_\Omega u^{\delta+r+1}dx
+4^{\delta+r}\mu|\Omega|,\\
\int_\Omega u^r dx&\le\frac{1}{4}\mu\int_\Omega u^{\delta+r+1}dx
+\Big(\frac{4}{\mu}\Big)^\frac{r}{\delta+r}|\Omega|,
\end{align*}
and
\begin{align} \nonumber
\int_\Omega u^m\phi(u)\nabla v\cdot \nabla u^r dx
&\le r\int_\Omega u^{m+r-1}|\nabla v\cdot\nabla u|dx
\\ \nonumber
&\le\frac{mr}{4}\int_\Omega(u+\varepsilon)^{m-1}u^{r-1}|\nabla u|^2dx
+\frac{r}{m}\int_\Omega u^{m+r}|\nabla v|^2dx
\\\label{eq-znablavMoser}
&\le\frac{1}{4}\int_\Omega\nabla(u+\varepsilon)^m\cdot\nabla u^r dx
+\frac{r}{m}\|\nabla v\|_{L^\infty(\Omega)}^2\int_\Omega u^{m+r}dx,
\end{align}
where according to Lemma \ref{le-nablavinfty}
$\|\nabla v\|_{L^\infty(\Omega)}$ is uniformly bounded.
Now we apply the Gagliardo-Nirenberg inequality Lemma \ref{le-GN} to obtain
\begin{align*}
\int_\Omega u^{m+r}dx&=\|u^\frac{m+r}{2}\|_{L^2(\Omega)}^2\\
&\le C_0\big(\|\nabla u^\frac{m+r}{2}\|_{L^2(\Omega)}^{2a}%
\|u^\frac{m+r}{2}\|_{L^1(\Omega)}^{2(1-a)}
+\|u^\frac{m+r}{2}\|_{L^1(\Omega)}^2\big),
\end{align*}
where $a=n/(n+2)\in(0,1)$
and $C_0$ is the constant in the Gagliardo-Nirenberg inequality
which is independent of $r$.
Therefore, we have
\begin{align} \nonumber
\frac{r}{m}\|\nabla v\|_{L^\infty(\Omega)}^2%
\int_\Omega u^{m+r}dx
&\le \frac{r}{m}\|\nabla v\|_{L^\infty(\Omega)}^2C_0
\big(\|\nabla u^\frac{m+r}{2}\|_{L^2(\Omega)}^{2a}%
\|u^\frac{m+r}{2}\|_{L^1(\Omega)}^{2(1-a)}
+\|u^\frac{m+r}{2}\|_{L^1(\Omega)}^2\big)\\ \nonumber
&\le\frac{mr}{(m+r)^2}\|\nabla u^\frac{m+r}{2}\|_{L^2(\Omega)}^2
+\Big(\frac{r}{m}\|\nabla v\|_{L^\infty(\Omega)}^2C_0\Big)^\frac{1}{1-a}%
\Big(\frac{(m+r)^2}{mr}\Big)^\frac{a}{1-a}
\|u^\frac{m+r}{2}\|_{L^1(\Omega)}^2 \\ \nonumber
&\qquad\qquad+\frac{r}{m}\|\nabla v\|_{L^\infty(\Omega)}^2C_0%
\|u^\frac{m+r}{2}\|_{L^1(\Omega)}^2\\ \label{eq-zumrmoser}
&\le\frac{1}{4}\int_\Omega\nabla(u+\varepsilon)^m\cdot\nabla u^r dx
+C_1(r)\|u^\frac{m+r}{2}\|_{L^1(\Omega)}^2,
\end{align}
where
$$C_1(r)=\Big(\frac{r}{m}\|\nabla v\|_{L^\infty(\Omega)}^2C_0\Big)^\frac{1}{1-a}%
\Big(\frac{(m+r)^2}{mr}\Big)^\frac{a}{1-a}
+\frac{r}{m}\|\nabla v\|_{L^\infty(\Omega)}^2C_0.$$
Inserting the above estimates \eqref{eq-znablavMoser}, \eqref{eq-zumrmoser} into
\eqref{eq-zur} yields
\begin{align} \nonumber
\frac{d}{dt}\int_\Omega u^{r+1}dx&+
\int_\Omega u^{r+1}dx \\ \nonumber
&\le C_1(r)(r+1)\|u^\frac{m+r}{2}\|_{L^1(\Omega)}^2
+(r+1)(4^{\delta+r}\mu|\Omega|+\Big(\frac{4}{\mu}\Big)^\frac{r}{\delta+r}|\Omega|)
\\ \nonumber
&\qquad\qquad+\int_\Omega u^{r+1}dx-\frac{1}{2}\mu\int_\Omega u^{\delta+r+1}dx
\\ \label{eq-zuLinfty}
&\le C_1(r)(r+1)\|u^\frac{m+r}{2}\|_{L^1(\Omega)}^2
+C_2(r),
\end{align}
where
$$C_2(r)=(r+1)\Big(4^{\delta+r}\mu|\Omega|
+\Big(\frac{4}{\mu}\Big)^\frac{r}{\delta+r}|\Omega|\Big)
+\Big(\frac{2}{\mu}\Big)^\frac{r+1}{\delta}|\Omega|.$$

Now we use the following Moser-type iteration.
Let $r=r_j$, with $r_j=2^j+m-2$ for $j\in\mathbb N^+$,
that is, $r_1=m$ and
$$r_{j-1}+1=\frac{r_j+m}{2}, \quad j\in \mathbb N^+.$$
We can invoke Lemma \ref{le-uLp} to find $C_0$ such that
$$\sup_{t\in(0,T_\text{max})}\|u\|_{L^{r_1+1}(\Omega)}\le C_0.$$
From \eqref{eq-zuLinfty} and an ODE comparison, we have
\begin{equation} \label{eq-zuLpj}
\sup_{t\in(0,T_\text{max})}\|u\|^{r_j+1}_{L^{r_j+1}(\Omega)}
\le\max\left\{\int_\Omega (u_0+1)^{r_j+1}dx,
C_1(r_j)(r_j+1)\cdot\sup_{t\in(0,T_\text{max})}
\|u\|_{L^{r_{j-1}+1}(\Omega)}^{2(r_{j-1}+1)}
+C_2(r_j)\right\}.
\end{equation}
A simple analysis shows that
$C_1(r)(r+1)\le a_1r^{b_1}$ and $C_2(r)\le a_2b_2^{r}$
for some positive constants $a_1,a_2$ and $b_1,b_2$
that all are greater than $1$ and independent of $r$.
Therefore, we can rewrite the above inequality \eqref{eq-zuLpj} into
\begin{equation} \label{eq-zuLpjnew}
\sup_{t\in(0,T_\text{max})}\|u\|^{r_j+1}_{L^{r_j+1}(\Omega)}
\le\max\left\{\int_\Omega (u_0+1)^{r_j+1}dx,
a_1r_j^{b_1}\cdot\sup_{t\in(0,T_\text{max})}
\|u\|_{L^{r_{j-1}+1}(\Omega)}^{2(r_{j-1}+1)}
+a_2b_2^{r_j}\right\}.
\end{equation}
Let
$$M_j=\max\Big\{\sup_{t\in(0,T_\text{max})}\int_\Omega u^{r_j+1}dx,1\Big\}.$$
Since boundedness of $u$ in $L^\infty(\Omega)$ is evident in the case when
$M_j\le\max\{\int_\Omega (u_0+1)^{r_j+1}dx,1\}$ for infinitely many $j\ge1$,
we may assume that $M_j\ge\max\{\int_\Omega (u_0+1)^{r_j+1}dx,1\}$ and thus,
according to \eqref{eq-zuLpjnew}, there holds
\begin{equation}\label{eq-zMj}
M_j\le a_1r_j^{b_1}M_{j-1}^2+a_2b_2^{r_j}.
\end{equation}
We note that if $M_{j-1}^2\le a_2b_2^{r_j}$ for infinitely many $j\ge1$, then
$$M_{j-1}^\frac{1}{r_{j-1}+1}\le (a_2b_2^{r_j})^\frac{1}{2(r_{j-1}+1)}
\le a_2^\frac{1}{r_j+m}b_2^\frac{r_j}{r_j+m}\le2b_2,$$
for $j$ sufficiently large, which shows the boundedness of $u$ in $L^\infty(\Omega)$.
Otherwise, $M_{j-1}^2\ge a_2b_2^{r_j}$ except for a finite number of $j\ge1$.
Thus, there exists a $j_0\ge1$ such that
$$M_{j-1}^2\ge a_2b_2^{r_j}, \qquad j\ge j_0.$$
Therefore, we can rewrite \eqref{eq-zMj} into
\begin{equation} \label{eq-zMjnew}
M_j\le2a_1r_j^{b_1}M_{j-1}^2\le D^jM_{j-1}^2
\end{equation}
for all $j\ge j_0$ with a constant $D$ independent of $j$,
whence upon enlarge $D$ if necessary we can achieve that
\eqref{eq-zMjnew} actually holds for all $j\ge1$.
By introduction, this yields
\begin{equation}\nonumber
M_j\le D^{\sum_{i=0}^{j-2}(j-i)\cdot 2^j}\cdot M_1^{2^{j-1}}
=D^{2^j+2^{j-1}-j-2}M_1^{2^{j-1}}\le D^{2^{j+1}}M_1^{2^{j-1}}
\end{equation}
for all $j\ge1$, and hence that
\begin{equation}\nonumber
M_j^{\frac{1}{r^j+1}}\le D^{\frac{2^{j+1}}{2^j+m-1}}
M_0^{\frac{2^{j-1}}{2^j+m-1}}
\le D^2 M_1,
\end{equation}
for all $j\ge1$.
This implies that $u$ indeed belongs to $L^\infty(\Omega\times(0,T_\text{max}))$.
$\hfill\Box$

Now we turn to the regularity estimates.

\begin{lemma} \label{le-wz}
Let $1\le n\le3$.
Then there exists a constant $C$ independent of $t$ and $\varepsilon$ such that
\begin{align*}
\sup_{t\in(0,T_\text{max})}(\|z_\varepsilon\|_{L^\infty(\Omega)}
+\|\nabla z_\varepsilon\|_{L^\infty(\Omega)}
+\|v_\varepsilon\|_{L^\infty(\Omega)}
+\|\nabla v_\varepsilon\|_{L^\infty(\Omega)})\le C.
\end{align*}
And the third solution component $w_\varepsilon$ fulfills
$$\|\nabla w_\varepsilon(\cdot,t)\|_{L^\infty(\Omega)}
\le2\|\nabla w_0\|_{L^\infty(\Omega)}
+(\|w_0\|_{L^\infty(\Omega)}+1)
\sup_{t\in(0,T_\text{max})}\|\nabla z_\varepsilon\|_{L^\infty(\Omega)}t,
\qquad t\in(0,T_\text{max}).$$
\end{lemma}
{\it \bfseries Proof.}
According to Lemma \ref{le-uLinfty}, Lemma \ref{le-zLq}, Lemma \ref{le-nablavLr},
we see that $\|u_\varepsilon\|_{L^\infty(\Omega}$,
$\|z_\varepsilon\|_{L^\infty(\Omega}$,
$\|\nabla v_\varepsilon\|_{L^\infty(\Omega}$
are uniformly bounded in $(0,T_\text{max})$.
The standard $L^p-L^q$ type estimates also shows the boundedness
of $\|v_\varepsilon\|_{L^\infty(\Omega}$ and
$\|\nabla z_\varepsilon\|_{L^\infty(\Omega}$.
We denote $v_\varepsilon,w_\varepsilon,z_\varepsilon$ by $v,w,z$
in this proof for the sake of simplicity.
Since both $w$ and $z$ are nonnegative
according to the third and fourth equation in \eqref{eq-model-regularized}
and the initial data, we have
\begin{align*}
w(x,t)&=w_{0\varepsilon}(x)e^{-\int_0^tz(x,\tau)d\tau}, \\
\nabla w(x,t)&=\nabla w_{0\varepsilon}(x)e^{-\int_0^tz(x,\tau)d\tau}
-w_{0\varepsilon}(x)e^{-\int_0^tz(x,\tau)d\tau}\int_0^t\nabla z(x,\tau)d\tau.
\end{align*}
Therefore,
\begin{align*}
|\nabla w(x,t)|&\le|\nabla w_{0\varepsilon}(x,t)|
+w_{0\varepsilon}(x)\sup_{t\in(0,T_\text{max})}\|\nabla z\|_{L^\infty(\Omega)}t\\
&\le2\|\nabla w_0\|_{L^\infty(\Omega)}
+(\|w_0\|_{L^\infty(\Omega)}+1)
\sup_{t\in(0,T_\text{max})}\|\nabla z_\varepsilon\|_{L^\infty(\Omega)}t.
\end{align*}
This completes the proof.
$\hfill\Box$

\begin{lemma} \label{le-Deltav}
There exists a constant $C>0$ independent of $\varepsilon$ and $T$, such that
$$\int_0^T\int_\Omega|\Delta v_\varepsilon|^2dxdt
\le C(1+T^2), \quad T\in(0,T_\text{max}).$$
\end{lemma}
{\it \bfseries Proof.}
We denote $v_\varepsilon,w_\varepsilon,z_\varepsilon$ by $v,w,z$
in this proof for the sake of simplicity.
Multiplying the second equation in \eqref{eq-model-regularized} by $-\Delta v$
and integrating over $\Omega$ yields
$$\int_\Omega\pd{}{t}|\nabla v|^2dx
+\int_\Omega|\Delta v|^2dx
=\int_\Omega\nabla v\cdot\nabla(wz)dx
\le C\Big(\int_\Omega|\nabla w|dx+1\Big)
\le C(1+t),$$
since $\nabla v$, $z$ and $\nabla z$ are uniformly bounded in $L^\infty(\Omega)$
according to Lemma \ref{le-wz}.
Integrating over $(0,T)$, we complete this proof.
$\hfill\Box$

\begin{lemma} \label{le-nablaum}
There exists a constant $C>0$ independent of $\varepsilon$ and $T$, such that
$$\int_0^T\int_\Omega |\nabla u_\varepsilon^m|^2dxdt\le C(1+T),
\quad T\in(0,T_\text{max}).$$
\end{lemma}
{\it \bfseries Proof.}
We denote $u_\varepsilon,v_\varepsilon$
by $u,v$ in this proof for the sake of simplicity.
We test the first equation in \eqref{eq-model-regularized}
by $(u+\varepsilon)^m$ and get
\begin{align} \nonumber
\frac{1}{m+1}\frac{d}{dt}\int_\Omega(u+\varepsilon)^{m+1}dx
&+\int_\Omega|\nabla(u+\varepsilon)^m|^2dx
\\ \nonumber
\le\int_\Omega u^m&\phi(u)\nabla v\cdot\nabla(u+\varepsilon)^mdx
+\mu\int_\Omega u^{\delta}(u+\varepsilon)^mdx
\\ \label{eq-znablaum}
&-\mu\int_\Omega u^{\delta+1}(u+\varepsilon)^mdx+\int_\Omega(u+\varepsilon)^mdx.
\end{align}
According to Lemma \ref{le-nablavinfty} and Lemma \ref{le-uLinfty},
$\nabla v$ and $u$ are uniformly bounded.
Thus,
\begin{align}\nonumber
\int_\Omega u^m\phi(u)\nabla v\cdot\nabla(u+\varepsilon)^mdx
\le\frac{1}{2}\int_\Omega|\nabla(u+\varepsilon)^m|^2dx
+C_1,
\end{align}
where $C_1$ is a constant independent of $t$ and $\varepsilon$.
Integrating \eqref{eq-znablaum} on $(0,T)$ yields
\begin{align} \label{eq-znablaumpe}
\int_\Omega (u+\varepsilon)^{m+1}dx
+\int_0^T\int_\Omega |\nabla(u+\varepsilon)^m|^2dx
\le\int_\Omega (u_{0\varepsilon}+\varepsilon)^{m+1}dx+CT.
\end{align}
We note that
$$|\nabla u^m|=mu^{m-1}|\nabla u|
\le m(u+\varepsilon)^{m-1}|\nabla (u+\varepsilon)|
=|\nabla(u+\varepsilon)^m|.$$
This completes the proof.
$\hfill\Box$

\begin{lemma} \label{le-umt}
There exists a constant $C>0$ independent of $\varepsilon$ and $T$, such that
$$\int_0^T\int_\Omega\Big|\Big(u_\varepsilon^\frac{m+1}{2}\Big)_t\Big|^2dxdt
+\int_\Omega\Big|\nabla u_\varepsilon^m\Big|^2dx
\le C(1+T^2),
\quad T\in(0,T_\text{max}).$$
Moreover,
$$\int_0^T\int_\Omega\Big|(u_\varepsilon^m)_t\Big|^2dxdt
\le\frac{4m^2}{(m+1)^2}\|u_\varepsilon\|_{L^\infty(\Omega)}^{m-1}%
\int_0^T\int_\Omega\Big|\Big(u_\varepsilon^\frac{m+1}{2}\Big)_t\Big|^2dxdt
\le C(1+T^2),
\quad T\in(0,T_\text{max}).$$
\end{lemma}
{\it \bfseries Proof.}
We denote $u_\varepsilon,v_\varepsilon$ by $u,v$
in this proof for the sake of simplicity.
We multiply the first equation in \eqref{eq-model-regularized}
by $[(u+\varepsilon)^m]_t$ and then we have
\begin{align} \nonumber
\int_\Omega m(u+\varepsilon)^{m-1}|u_t|^2dx
&+\int_\Omega\nabla(u+\varepsilon)^m\cdot\nabla[(u+\varepsilon)^m]_tdx
\\ \nonumber
\le\int_\Omega u^m&\phi(u)\nabla v\cdot\nabla[(u+\varepsilon)^m]_tdx
+\mu\int_\Omega u^{\delta}[(u+\varepsilon)^m]_tdx
\\ \label{eq-zumt}
&-\mu\int_\Omega u^{\delta+1}[(u+\varepsilon)^m]_tds
+\int_\Omega\big|[(u+\varepsilon)^m]_t\big|dx.
\end{align}
We note that $\|u\|_{L^\infty(\Omega)}$ is uniformly bounded and then
\begin{align*}
\int_\Omega\mu u^{\delta}[(u+\varepsilon)^m]_tdx
&=\int_\Omega m\mu u^\delta(u+\varepsilon)^{m-1}u_tdx
\le\frac{1}{5}\int_\Omega m(u+\varepsilon)^{m-1}|u_t|^2dx
+C_1, \\
\int_\Omega-\mu u^{\delta+1}[(u+\varepsilon)^m]_tdx
&=-\int_\Omega m\mu u^{\delta+1}(u+\varepsilon)^{m-1}u_tdx
\le\frac{1}{5}\int_\Omega m(u+\varepsilon)^{m-1}|u_t|^2dx
+C_2, \\
\int_\Omega\big|[(u+\varepsilon)^m]_t\big|dx
&=\int_\Omega m(u+\varepsilon)^{m-1}u_tdx
\le\frac{1}{5}\int_\Omega m(u+\varepsilon)^{m-1}|u_t|^2dx
+C_3,
\end{align*}
where $C_1,C_2,C_3$ are constants independent of $t$ and $\varepsilon$.
We also have
\begin{align*}
\int_\Omega m(u+\varepsilon)^{m-1}|u_t|^2dx
=\frac{4m}{(m+1)^2}\int_\Omega%
\Big|\Big((u+\varepsilon)^\frac{m+1}{2}\Big)_t\Big|^2dx,
\end{align*}
and
\begin{align*}
\int_\Omega\nabla(u+\varepsilon)^m\cdot\nabla[(u+\varepsilon)^m]_tdx
=\frac{1}{2}\pd{}{t}\int_\Omega\Big|\nabla(u+\varepsilon)^m\Big|^2dx.
\end{align*}
There holds
\begin{align*}
&\int_\Omega u^m\phi(u)\nabla v\cdot\nabla[(u+\varepsilon)^m]_tdx
=-\int_\Omega[(u+\varepsilon)^m]_t\nabla\cdot(u^m\phi(u)\nabla v)dx \\
=&-\int_\Omega m(u+\varepsilon)^{m-1}u_t\cdot%
(mu^{m-1}\phi(u)\nabla u\cdot\nabla v
+u^m\phi'(u)\nabla u\cdot\nabla v
+u^m\phi(u)\Delta v)dx\\
\le&\frac{1}{5}\int_\Omega m(u+\varepsilon)^{m-1}|u_t|^2dx
+C_4\int_\Omega(u+\varepsilon)^{2(m-1)}|\nabla u|^2dx
+C_5\int_\Omega|\Delta v|^2dx\\
\le&\frac{1}{5}\int_\Omega m(u+\varepsilon)^{m-1}|u_t|^2dx
+C_4\int_\Omega|\nabla(u+\varepsilon)^m|^2dx
+C_5\int_\Omega|\Delta v|^2dx,
\end{align*}
where $C_4$ and $C_5$ are constants independent of $t$ and $\varepsilon$,
since the uniform boundedness of $\|\nabla v\|_{L^\infty(\Omega)}$.
Inserting the above inequalities into \eqref{eq-zumt},
and noticing the inequality \eqref{eq-znablaumpe}
in the proof of Lemma \ref{le-nablaum},
we find a constant $C$ independent of $t$ and $\varepsilon$ such that
\begin{equation*}
\int_0^T\int_\Omega\Big|\Big((u+\varepsilon)^\frac{m+1}{2}\Big)_t\Big|^2dxdt
+\int_\Omega\Big|\nabla(u+\varepsilon)^m\Big|^2dx
\le\int_\Omega\Big|\nabla(u_{0\varepsilon}+\varepsilon)^m\Big|^2dx
+C(1+T^2)\le C(1+T^2).
\end{equation*}
Clearly, we have
$$\Big|\Big(u^\frac{m+1}{2}\Big)_t\Big|^2
=\frac{(m+1)^2}{4}u^{m-1}|u_t|^2
\le\frac{(m+1)^2}{4}(u+\varepsilon)^{m-1}|u_t|^2
=\Big|\Big((u+\varepsilon)^\frac{m+1}{2}\Big)_t\Big|^2,$$
and
$$|(u^m)_t|^2\le\frac{4m^2}{(m+1)^2}\|u_\varepsilon\|_{L^\infty(\Omega)}^{m-1}%
\Big|\Big(u^\frac{m+1}{2}\Big)_t\Big|^2
\le\frac{4m^2}{(m+1)^2}\|u_\varepsilon\|_{L^\infty(\Omega)}^{m-1}%
\Big|\Big((u+\varepsilon)^\frac{m+1}{2}\Big)_t\Big|^2.$$
The proof is completed.
$\hfill\Box$

\noindent
{\it \bfseries Proof of Theorem \ref{th-existence}.}
According to the estimates, for any $\varepsilon$, the approximation solution
$(u_\varepsilon,v_\varepsilon,w_\varepsilon,z_\varepsilon)$ exists globally.
The regularity estimates of $v_\varepsilon,w_\varepsilon,z_\varepsilon$ are trivial.
For any $T\in(0,\infty)$, we see that
$u_\varepsilon^m\in L^\infty(Q_T)$, $\nabla u_\varepsilon^m\in L^2(Q_T)$,
and $\partial{u_\varepsilon^m}/{\partial t}\in L^2(Q_T)$,
Thus, there exists a function $\widetilde u\in W^{1,2}(Q_T)$, such that
$u_\varepsilon^m$ weakly in $W^{1,2}(Q_T)$
and strongly in $L^2(Q_T)$ converges to $\widetilde u$.
We denote $u=\widetilde u^{1/m}$ since $\widetilde u\ge0$.
Thus, $u_\varepsilon^m$ converges almost everywhere to $u^m$, and
$u_\varepsilon$ converges almost everywhere to $u$.
We can verify the integral identities in the definition of weak solutions.
By taking a sequence of $T\in(0,\infty)$ and the diagonal subsequence procedure,
we can find the existence of a global weak solution.
$\hfill\Box$

Now we show the regularity of the globally bounded weak solution.

\begin{lemma} \label{le-bounded}
Let $(u,v,w,z)$ be a globally bounded weak solution of \eqref{eq-model}
such that $\|u(\cdot,t)\|_{L^\gamma(\Omega)}$ is uniformly bounded
with $\gamma=\max\{1,n/3\}$.
Then there exists a constant $C$ such that
$$
\sup_{t\in\mathbb R^+}\left\{\|u\|_{L^\infty(\Omega)}+\|v\|_{W^{1,\infty}(\Omega)}
+\|w\|_{L^{\infty}(\Omega)}+\|z\|_{W^{1,\infty}(\Omega)}\right\}
\le C.
$$
\end{lemma}
{\it\bfseries Proof.}
Since $\|u(\cdot,t)\|_{L^\frac{n}{3}(\Omega)}$ is uniformly bounded,
for any $r\ge1$ we can apply Lemma \ref{le-zLq} and \ref{le-nablavLr} to
find a constant $C(r)$ independent of $t$ such that
$\|\nabla v(\cdot,t)\|_{L^r(\Omega)}\le C(r)$ for all $t>0$.
The estimates in the proof of Theorem \ref{th-existence} in section 3
can be carried on to complete this proof here.
$\hfill\Box$

In Lemma \ref{le-wz}, we have proved
$\|\nabla w(\cdot,t)\|_{L^\infty(\Omega)}\le C(1+t)$
(same as $\nabla w_\varepsilon$) for some constant $C>0$.
However that is an estimate depending on time $t$.
Employing the method in the proof of Lemma \ref{le-nablaw} in next Section and iteration technique,
we can prove the following uniform estimate.

\begin{lemma} \label{le-nablawp}
Let $(u,v,w,z)$ be a globally bounded weak solution of \eqref{eq-model}.
Then for any $p\ge1$ there holds
$$\int_\Omega|\nabla w(\cdot,t)|^pdx
\le C(p), \qquad t>0,$$
for some constant $C(p)$ independent of time $t$.
\end{lemma}
{\it\bfseries Proof.}
This proof proceeds along the idea of the arguments of
Lemma 4.3 in \cite{WangLargeJDE} and Lemma 4.1 in \cite{TaoWinkler}.
Since
$$w(x,t)=w_0(x)e^{-\int_0^tz(x,s)ds},$$
and
$$\nabla w(x,t)=\nabla w_0(x)e^{-\int_0^tz(x,s)ds}
-w_0(x)e^{-\int_0^tz(x,s)ds}\int_0^t\nabla z(x,s)ds.$$
We see that
\begin{align*}
|\nabla w(x,t)|^2
\le2|\nabla w_0(x)|^2e^{-2\int_0^tz(x,s)ds}
+2|w_0(x)|^2e^{-2\int_0^tz(x,s)ds}\Big|\int_0^t\nabla z(x,s)ds\Big|^2.
\end{align*}
And thus
\begin{align*}
\int_\Omega|\nabla w(x,t)|^2dx
\le& C+C\int_\Omega e^{-2\int_0^tz(x,s)ds}\Big|\int_0^t\nabla z(x,s)ds\Big|^2dx\\
\le& C-\frac{C}{2}\int_\Omega\nabla e^{-2\int_0^tz(x,s)ds}%
\cdot\Big(\int_0^t\nabla z(x,s)ds\Big)dx\\
\le& C+\frac{C}{2}\int_\Omega e^{-2\int_0^tz(x,s)ds}%
\cdot\Big(\int_0^t\Delta z(x,s)ds\Big)dx\\
\le& C+\frac{C}{2}\int_\Omega e^{-2\int_0^tz(x,s)ds}%
\cdot\Big(\int_0^t(z_t+z-u)ds\Big)dx\\
\le& C+\frac{C}{2}\int_\Omega e^{-2\int_0^tz(x,s)ds}%
\cdot\Big(z(x,t)+\int_0^tz(x,s)ds\Big)dx\\
\le& C.
\end{align*}
Using the same method, we have
\begin{align*}
|\nabla w(x,t)|^4
\le2^3|\nabla w_0(x)|^4e^{-4\int_0^tz(x,s)ds}
+2^3|w_0(x)|^4e^{-4\int_0^tz(x,s)ds}\Big|\int_0^t\nabla z(x,s)ds\Big|^4,
\end{align*}
and
\begin{align*}
\int_\Omega|\nabla w(x,t)|^4dx
\le& C+C\int_\Omega e^{-4\int_0^tz(x,s)ds}\Big|\int_0^t\nabla z(x,s)ds\Big|^4dx\\
\le& C-\frac{C}{4}\int_\Omega\nabla e^{-4\int_0^tz(x,s)ds}%
\cdot\Big(\int_0^t\nabla z(x,s)ds\Big)^3dx\\
\le& C+\frac{3C}{4}\int_\Omega e^{-4\int_0^tz(x,s)ds}%
\cdot\Big(\int_0^t\nabla z(x,s)ds\Big)^2%
\cdot\Big(\int_0^t\Delta z(x,s)ds\Big)dx\\
\le& C+\frac{3C}{4}\int_\Omega e^{-4\int_0^tz(x,s)ds}%
\cdot\Big(\int_0^t\nabla z(x,s)ds\Big)^2%
\cdot\Big(\int_0^t(z_t+z-u)ds\Big)dx\\
\le& C+\frac{3C}{4}\int_\Omega e^{-4\int_0^tz(x,s)ds}%
\cdot\Big(\int_0^t\nabla z(x,s)ds\Big)^2%
\cdot\Big(z(x,t)+\int_0^tz(x,s)ds\Big)dx\\
\le& C+\frac{3C}{4}\int_\Omega e^{-2\int_0^tz(x,s)ds}%
\cdot\Big(\int_0^t\nabla z(x,s)ds\Big)^2dx\\
&\qquad\cdot\sup_{x\in\Omega}\Big[e^{-2\int_0^tz(x,s)ds}%
\cdot\Big(z(x,t)+\int_0^tz(x,s)ds\Big)\big]\\
\le&C,
\end{align*}
according to the proof of the previous estimate on $\|\nabla w\|_{L^2(\Omega)}$
and the boundedness of $\|z\|_{L^\infty(\Omega)}$.
Repeating this process for $\|\nabla w\|_{L^k(\Omega)}$ with $k=6,8,\dots$,
we complete this proof by iteration.
$\hfill\Box$

\noindent
{\it \bfseries Proof of Theorem \ref{th-regular}.}
Lemma \ref{le-nablawp} shows the uniform boundedness of
$\|\nabla w\|_{L^{n+2}(\Omega)}$.
According to the third equation of \eqref{eq-model},
we see that $\|w_t\|_{L^\infty(\Omega)}=\|wz\|_{L^\infty(\Omega)}\le C$.
Therefore, $w\in W^{1,n+2}(\Omega\times(t,t+1))$
and its norm is uniformly bounded for any $t>0$.
Sobolev embedding theorem implies the existence of
$\alpha\in(0,1)$ and $C>0$ such that
$$\|w\|_{C^\alpha(\overline\Omega\times[t,t+1])}\le C,\quad t>0.$$
Since $u$ is uniformly bounded,
the strong solution theory of parabolic equation
applied to the fourth equation in \eqref{eq-model} shows
$$\|z_t\|_{L^p(\Omega\times(t,t+1))}
+\|\Delta z\|_{L^p(\Omega\times(t,t+1))}\le C(p),\quad t>0,$$
for some constant $C(p)>0$.
Taking $p>1+n/2$, we see that for some $\alpha\in(0,1)$
$$\|z\|_{C^\alpha(\overline\Omega\times[t,t+1])}\le C,\quad t>0.$$
Thus,
$$\|wz\|_{C^\alpha(\overline\Omega\times[t,t+1])}\le C,\quad t>0.$$
This can also be deduced by
$$\|\nabla(wz)\|_{L^p(\Omega)}+\|(wz)_t\|_{L^p(\Omega\times(t,t+1))}\le C,\quad t>0,$$
with $p>n+1$.
Using bootstrap arguments involving the standard parabolic regularity theory,
we can verify that
$$\|v\|_{C^{2+\alpha,1+\alpha/2}(\overline\Omega\times[t,t+1])}
+\|z\|_{W^{2,1}_{p}(\Omega\times(t,t+1))}\le C(p).$$
The proof is completed.
$\hfill\Box$

\section{Propagating properties and large time behavier}

This section is devoted to the study of the propagating properties of the tumour cells and
the large time behavior of the weak solution
$(u,v,w,z)$ to the problem \eqref{eq-model}.
In contrast with the heat equation,
it is known that the porous medium equation has the
property of finite speed of propagation.
Therefore, the first component $u$ may not have positive minimum
for some time $t>0$.
We use the comparison principle together with two kinds of weak lower solutions,
one is decaying but its support is expanding with finite speed of propagation,
the other one is an increasing function of time $t$,
to overcome the difficulty of degenerate dispersion.

We first present the following comparison principle of
the first component.

\begin{lemma} \label{le-comparisonprinciple}
Let $T>0$ and the function space
$$E=\{u\in L^\infty(Q_T);u\ge0,\nabla u^m\in L^2((0,T);L^2(\Omega)),
u^{m-1}u_t\in L^2((0,T);L^2(\Omega))\},$$
$u_1,u_2\in E$, $\nabla v\in L^\infty(Q_T)$,
and $u_1$, $u_2$ satisfy the following differential inequalities
\begin{align*}
\begin{cases}
\displaystyle
\pd{u_1}{t}\ge\Delta u_1^m-\nabla\cdot(u_1^m\phi(u_1)\nabla v)
+\mu u_1^\delta(1-u_1), &\\[2mm]
\displaystyle
\pd{u_2}{t}\le\Delta u_2^m-\nabla\cdot(u_2^m\phi(u_2)\nabla v)
+\mu u_2^\delta(1-u_2),\qquad &x\in\Omega, t\in(0,T), \\[2mm]
\displaystyle
\pd{u_1}{\nu}\ge0\ge\pd{u_2}{\nu}, \qquad &x\in\partial\Omega, t\in(0,T),\\
u_1(x,0)\ge u_2(x,0)\ge0, \qquad &x\in\Omega,
\end{cases}
\end{align*}
in the sense that the following inequalities
\begin{align*}
\int_0^T\int_\Omega u_1\varphi_tdxdt&
+\int_\Omega u_{10}(x)\varphi(x,0) dx
\le\int_0^T\int_\Omega \nabla u_1^m\cdot\nabla\varphi dxdt\\
&\qquad-\int_0^T\int_\Omega u_1^m\phi(u_1)\nabla v\cdot\nabla\varphi dxdt
-\int_0^T\int_\Omega\mu u_1^\delta(1-u_1)\varphi dxdt, \\
\int_0^T\int_\Omega u_2\varphi_tdxdt&
+\int_\Omega u_{20}(x)\varphi(x,0) dx
\ge\int_0^T\int_\Omega \nabla u_2^m\cdot\nabla\varphi dxdt\\
&\qquad-\int_0^T\int_\Omega u_2^m\phi(u_2)\nabla v\cdot\nabla\varphi dxdt
-\int_0^T\int_\Omega\mu u_2^\delta(1-u_2)\varphi dxdt,
\end{align*}
hold for some fixed $u_{10},u_{20}\in L^2(\Omega)$
such that $u_{10}\ge u_{20}\ge0$ on $\Omega$ and all test functions
$0\le\varphi\in L^2((0,T);W^{1,2}(\Omega))\cap W^{1,2}((0,T);L^2(\Omega))$ with
$\varphi(x,T)=0$ on $\Omega$.
Then $u_1(x,t)\ge u_2(x,t)$ almost everywhere in $Q_T$.
\end{lemma}
{\it \bfseries Proof.}
The following inequality
\begin{align*}
\int_0^T&\int_\Omega(u_1-u_2)\varphi_tdxdt
\le\int_0^T\int_\Omega \nabla (u_1^m-u_2^m)\cdot\nabla\varphi dxdt\\
&\qquad-\int_0^T\int_\Omega (u_1^m\phi(u_1)-u_2^m\phi(u_2))
\nabla v\cdot\nabla\varphi dxdt
-\int_0^T\int_\Omega\mu (u_1^\delta(1-u_1)-u_2^\delta(1-u_2))\varphi dxdt,
\end{align*}
holds for all
$0\le\varphi\in L^2((0,T);W^{1,2}(\Omega))\cap W^{1,2}((0,T);L^2(\Omega))$ with
$\varphi(x,T)=0$.
Let
$$a(x,t)=
\begin{cases}
\displaystyle
\frac{u_1^m-u_2^m}{u_1-u_2}, \quad &u_1(x,t)\ne u_2(x,t),\\
mu_1^{m-1}, \quad &u_1(x,t)=u_2(x,t),
\end{cases}$$
$$b(x,t)=
\begin{cases}
\displaystyle
\frac{(u_1^m\phi(u_1)-u_2^m\phi(u_2))\nabla v}{u_1-u_2},
\quad &u_1(x,t)\ne u_2(x,t),\\
(mu_1^{m-1}\phi(u_1)+u_1^m\phi'(u_1))\nabla v, \quad &u_1(x,t)=u_2(x,t),
\end{cases}$$
and
$$c(x,t)=
\begin{cases}
\displaystyle
\frac{\mu (u_1^\delta(1-u_1)-u_2^\delta(1-u_2))}{u_1-u_2},
\quad &u_1(x,t)\ne u_2(x,t),\\
\mu\delta u_1^{\delta-1}-\mu(\delta+1)u_1^\delta,
\quad &u_1(x,t)=u_2(x,t).
\end{cases}$$
Since $\nabla v,u_1,u_2$ are bounded and $\phi$ is smooth enough, there exists a constant
$C>0$ such that $|b|\le C a$ and $|c|\le C$.
Henceforth, a generic positive constant (possibly changing from line to line)
is denoted by $C$.
However, $c$ is not bounded by $Ca$
and we have no estimate on $\nabla c$
since we only assume that $\delta\ge1$.
Then for all
$0\le\varphi\in L^2((0,T);W^{1,2}(\Omega))\cap W^{1,2}((0,T);L^2(\Omega))$ with
$\varphi(x,T)=0$ on $\Omega$
and $\pd{\varphi}{\nu}=0$ on $\partial\Omega\times(0,T)$, there holds
\begin{align*}
\int_0^T&\int_\Omega(u_1-u_2)
(\varphi_t+a(x,t)\Delta\varphi+b(x,t)\cdot\nabla\varphi+c(x,t)\varphi)
dxdt\le0.
\end{align*}
We employ the standard duality proof method
or the approximate Hohmgren's approach to complete this proof
(see Theorem 6.5 in \cite{PME}, Chapter 1.3 and 3.2 in \cite{NDE}).
For any smooth function $\psi(x,t)\ge0$, we solve the inverse-time problem
\begin{equation} \label{eq-zdual}
\begin{cases}
\varphi_t+(\kappa+a_\varepsilon(x,t))\Delta\varphi
+b(x,t)\cdot\nabla\varphi+c_\theta(x,t)\varphi+\psi=0,
\quad &(x,t)\in Q_T, \\
\displaystyle
\pd{\varphi}{\nu}=0, \quad &(x,t)\in\partial\Omega\times(0,T),\\
\varphi(x,T)=0, \quad &x\in\Omega,
\end{cases}
\end{equation}
where $\kappa>0$, $\theta>0$, $a_\varepsilon$ is a smooth approximation of $a$,
$a_\varepsilon\ge a$, and
$$c_\theta(x,t)=
\begin{cases}
\displaystyle
\frac{\mu (u_1^\delta(1-u_1)-u_2^\delta(1-u_2))}{u_1-u_2},
\quad &|u_1(x,t)-u_2(x,t)|\ge\theta,\\
0, \quad &|u_1(x,t)-u_2(x,t)|<\theta.
\end{cases}$$
This definition of $c_\theta$ allows us to find a constant $C(\theta)$ such that
$$\frac{c_\theta^2}{a}\le C(\theta).$$
We may also need to replace $b(x,t)$ and $c_\theta(x,t)$ by
their smooth approximation functions
$b_\varepsilon(x,t)$ and $c_{\theta,\varepsilon}(x,t)$ respectively
in \eqref{eq-zdual}.
For the sake of simplicity we omit this procedure.
Here we note that \eqref{eq-zdual} is a standard parabolic problem
as the initial data is imposed at the end time $t=T$.
Therefore, it has a smooth solution $\varphi\ge0$.
Maximum principle shows the boundedness of $\varphi$.
Then we get the estimate
\begin{align*}
\iint_{Q_T}(u_1-u_2)\psi dxdt\ge&
-\iint_{Q_T}|u_1-u_2||a-a_\varepsilon||\Delta\varphi|dxdt\\
&-\kappa\iint_{Q_T}|u_1-u_2||\Delta\varphi|dxdt
-\iint_{Q_T}|u_1-u_2||c-c_\theta|\varphi dxdt\\
=:&-I_1-I_2-I_3.
\end{align*}
Now we need the a priori estimate on $a_\varepsilon|\Delta\varphi|^2$.
We can assume that $T$ is appropriately small,
otherwise we can prove step by step on each time interval.
We multiply the equation \eqref{eq-zdual} by $\eta(t)\Delta\varphi$
where $1/2\le\eta(t)\le1$ is a smooth function with $\eta'(t)\ge M>0$ for $t\in(0,T)$.
Since $T$ is small, we can choose $M$ appropriately large.
Integrating over $Q_T$ yields
\begin{align*}
&\iint_{Q_T}\varphi_t\eta\Delta\varphi dxdt
+\iint_{Q_T}\eta (\kappa+a_\varepsilon)(\Delta\varphi)^2dxdt\\
\le& \iint_{Q_T}\eta |b||\nabla\varphi||\Delta\varphi|dxdt
+\iint_{Q_T}\eta c_\theta\varphi\Delta\varphi dxdt
+\iint_{Q_T}\eta\psi\Delta\varphi dxdt\\
\le&\iint_{Q_T}\eta Ca|\nabla\varphi||\Delta\varphi|dxdt
+\frac{1}{4}\iint_{Q_T}\eta (\kappa+a_\varepsilon)(\Delta\varphi)^2dxdt\\
&\qquad+\iint_{Q_T}\frac{\eta c_\theta^2\varphi^2}{\kappa+a_\varepsilon}dxdt
+\iint_{Q_T}\eta|\nabla\psi||\nabla\varphi|dxdt\\
\le&\frac{1}{2}\iint_{Q_T}\eta (\kappa+a_\varepsilon)(\Delta\varphi)^2dxdt
+\iint_{Q_T}\frac{\eta C^2a^2|\nabla\varphi|^2}{\kappa+a_\varepsilon}dxdt
+\iint_{Q_T}\eta C(\theta)\varphi^2dxdt\\
&\qquad+\iint_{Q_T}\eta|\nabla\psi|^2dxdt
+\iint_{Q_T}\eta|\nabla\varphi|^2dxdt.
\end{align*}
Using $\varphi(x,T)=0$, we have
\begin{align*}
\iint_{Q_T}\varphi_t\eta\Delta\varphi dxdt
&=-\iint_{Q_T}\eta\nabla\varphi\cdot\nabla\varphi_t dxdt
=-\frac{1}{2}\iint_{Q_T}\eta\pd{}{t}|\nabla\varphi|^2 dxdt\\
&\ge\frac{1}{2}\iint_{Q_T}\eta'(t)|\nabla\varphi|^2 dxdt
\ge\frac{M}{2}\iint_{Q_T}|\nabla\varphi|^2 dxdt.
\end{align*}
Therefore,
\begin{align} \label{eq-zaDeltaphi}
\iint_{Q_T}|\nabla\varphi|^2 dxdt
+\iint_{Q_T}(\kappa+a_\varepsilon)(\Delta\varphi)^2dxdt
\le C(\theta).
\end{align}
It follows that
\begin{align*}
I_1=&\iint_{Q_T}|u_1-u_2||a-a_\varepsilon||\Delta\varphi|dxdt\\
\le&\Big(\iint_{Q_T}(\kappa+a_\varepsilon)|\Delta\varphi|^2dxdt\Big)^\frac{1}{2}
\cdot\Big(\iint_{Q_T}\frac{|a-a_\varepsilon|^2}{\kappa+a_\varepsilon}
|u_1-u_2|^2dxdt\Big)^\frac{1}{2}\\
\le&C(\theta)\Big(\iint_{Q_T}\frac{|a-a_\varepsilon|^2}{\kappa+a_\varepsilon}
dxdt\Big)^\frac{1}{2}\\
\le&\frac{C(\theta)}{\kappa^\frac{1}{2}}
\Big(\iint_{Q_T}|a-a_\varepsilon|^2dxdt\Big)^\frac{1}{2},
\end{align*}
which converges to zero if we let $\varepsilon\to0$.
For any fixed $\gamma>0$, denote
$$F_\gamma=\{(x,t)\in Q_T;|u_1-u_2|\ge\gamma\},$$
and
$$G_\gamma=\{(x,t)\in Q_T;|u_1-u_2|<\gamma\}.$$
Then there exists a constant $C(\gamma)$ such that
$a(x,t)\ge C(\gamma)$ on $F_\gamma$ and
\begin{align*}
I_2&=\kappa\iint_{Q_T}|u_1-u_2||\Delta\varphi|dxdt\\
&\le\kappa\iint_{G_\gamma}|u_1-u_2||\Delta\varphi|dxdt
+\kappa\iint_{F_\gamma}|u_1-u_2||\Delta\varphi|dxdt\\
&\le\gamma\iint_{G_\gamma}\kappa|\Delta\varphi|dxdt
+\frac{C\kappa}{C(\gamma)^\frac{1}{2}}%
\iint_{F_\gamma}a^\frac{1}{2}|\Delta\varphi|dxdt\\
&\le C\gamma\Big(\iint_{Q_T}\kappa|\Delta\varphi|^2dxdt\Big)^\frac{1}{2}
+\frac{C\kappa}{C(\gamma)^\frac{1}{2}}%
\Big(\iint_{Q_T}a|\Delta\varphi|^2dxdt\Big)^\frac{1}{2}\\
&\le \gamma C(\theta)+\frac{\kappa C(\theta)}{C(\gamma)^\frac{1}{2}},
\end{align*}
which converges to zero if we first let $\kappa\to0$ and then let $\gamma\to0$.
We also have
\begin{align*}
I_3=\iint_{Q_T}|u_1-u_2||c-c_\theta|\varphi dxdt
\le C\Big(\iint_{Q_T}|c-c_\theta|^2dxdt\Big)^\frac{1}{2},
\end{align*}
which converges to zero if we let $\theta\to0$.
Now we conclude that
\begin{align*}
\iint_{Q_T}(u_1-u_2)\psi dxdt\ge0
\end{align*}
for any given $\psi\ge0$ and then $u_1\ge u_2$ almost everywhere on $Q_T$.
$\hfill\Box$

Here we recall some lemmas about the asymptotic behavior of solutions to
evolutionary equations.

\begin{lemma}[\cite{Ito}] \label{le-v}
Let $(u,v,w,z)$ be a global solution of \eqref{eq-model}.
Then there exists a constant $L\ge0$ such that
$$\|v(\cdot,t)-L\|_{W^{1,\infty}(\Omega)}\to0,
\qquad\text{as~}t\to\infty.$$
In particular,
$$\|\nabla v(\cdot,t)\|_{L^{\infty}(\Omega)}\to0,
\qquad\text{as~}t\to\infty.$$
\end{lemma}

\begin{lemma}[\cite{WangLargeJDE} Lemma 4.1] \label{le-intv}
If $z$ is a global classical solution of
\begin{align*}
\begin{cases}
z_t=\Delta z-z+u, \quad &x\in\Omega, ~t>0,\\
\displaystyle
\pd{z}{\nu}=0, \quad &x\in\partial\Omega, ~t>0,\\
z(x,0)=z_0(x), \quad &x\in\Omega,
\end{cases}
\end{align*}
where $u(x,t)\ge0$ is given.
Then there exist constants $C_1$ and $C_2>0$
only depend on $\text{diam}\Omega$ and $\sup_{\tau<t}\|u\|_{L^1(\Omega)}$
respectively, such that
$$\int_0^tz(x,s)ds\ge
C_1\int_0^t\int_\Omega u(y,s)dyds-C_2, \quad x\in\Omega,~t>0.$$
\end{lemma}

\begin{lemma}[\cite{WangLargeJDE} Lemma 4.3, \cite{TaoWinkler} Lemma 4.1] \label{le-nablaw}
If $(w,z)$ is a global solution of
\begin{align*}
\begin{cases}
w_t=-wz, \\
z_t=\Delta z-z+u, \quad &x\in\Omega, ~t>0,\\[1mm]
\displaystyle
\pd{z}{\nu}=0, \quad &x\in\partial\Omega, ~t>0,\\
w(x,0)=w_0(x), \\
z(x,0)=z_0(x), \quad &x\in\Omega,
\end{cases}
\end{align*}
with $u\ge0$ on $\Omega\times\mathbb R^+$
and $\pd{w_0}{\nu}=0$ on $\partial\Omega$, then
$$\int_\Omega|\nabla w(\cdot,t)|^2dx
\le2\int_\Omega|\nabla w_0|^2dx
+\frac{|\Omega|}{2e}\|w_0\|_{L^\infty(\Omega)}^2
+\|w_0\|_{L^\infty(\Omega)}^2\int_\Omega z(\cdot,t)dx$$
for all $t>0$.
\end{lemma}

Now we construct a self similar weak lower solution with expanding support.

\begin{lemma} \label{le-lower}
Let $(u,v,w,z)$ be a globally bounded weak solution of \eqref{eq-model}
with the first component initial data $u_0\ge0$, $u_0\not\equiv0$
and $1\le\delta<m$, $\Omega$ is convex.
Define a function
$$g(x,t)=\varepsilon(1+t)^{-\kappa}
\Big[\Big(\eta-\frac{|x-x_0|^2}{(1+t)^\beta}\Big)_+\Big]^d,
\quad x\in\Omega,~t\ge0,$$
where $d=1/(m-1)$, $\beta\in(0,1/2)$ is sufficiently small, $\kappa=(1-\beta)/(m-1)$,
$x_0\in\Omega$ such that $\inf_{x\in B_{r}(x_0)}u_0(x)>0$ for some $r>0$,
$\varepsilon\in(0,1/2)$, $\eta>0$.
Then by appropriately selecting $\beta$, $\varepsilon$ and $\eta$,
the function $g(x,t)$ can be
a weak lower solution of the first equation in \eqref{eq-model},
that is,
\begin{align*}
\begin{cases}
\displaystyle
\pd{g}{t}\le\Delta g^m-\nabla\cdot(g^m\phi(g)\nabla v)
+\mu g^\delta(1-g),\qquad &x\in\Omega, t\in(0,T), \\[3mm]
\displaystyle
\pd{g}{\nu}\le0, \qquad &x\in\partial\Omega, t\in(0,T),\\[2mm]
0\le g(x,0)\le u_0(x), \qquad &x\in\Omega,
\end{cases}
\end{align*}
in the sense that the following inequality
\begin{align*}
\int_0^T\int_\Omega g\varphi_tdxdt&
+\int_\Omega g(x,0)\varphi(x,0) dx
\ge\int_0^T\int_\Omega \nabla g^m\cdot\nabla\varphi dxdt\\
&\qquad-\int_0^T\int_\Omega g^m\phi(g)\nabla v\cdot\nabla\varphi dxdt
-\int_0^T\int_\Omega\mu g^\delta(1-g)\varphi dxdt,
\end{align*}
holds for any $T>0$ and all test functions
$0\le\varphi\in L^2((0,T);W^{1,2}(\Omega))\cap W^{1,2}((0,T);L^2(\Omega))$ with
$\varphi(x,T)=0$ on $\Omega$,
and $0\le g(x,0)\le u_0(x)$ on $\Omega$.
Therefore, $u(x,t)\ge g(x,t)$ and there exist $t_0>0$ and $\varepsilon_0\ge0$
such that $u(x,t)\ge\varepsilon_0$ for all $x\in\Omega$ and $t\ge t_0$.
\end{lemma}
{\it\bfseries Proof.}
For simplicity, we let
$$h(x,t)=\Big(\eta-\frac{|x-x_0|^2}{(1+t)^\beta}\Big)_+, \quad x\in\Omega, ~t\ge0,$$
and
$$A(t)=\Big\{x\in\Omega;\frac{|x-x_0|^2}{(1+t)^\beta}<\eta\Big\},
\quad t\ge0.$$
Since $u_0\ge0$, $u_0\not\equiv0$ and $u_0\in C(\overline\Omega)$,
we see that there exists
$x_0\in\Omega$ such that $u_0(x)\ge\varepsilon_1$ on $B_{r}(x_0)$
for some $r>0$ and $\varepsilon_1>0$.
Without loss of generality, we may assume that $B_{r}(x_0)\subset\Omega$,
$x_0=0$ and $\varepsilon_1\le1/2$.
Straightforward computation shows that
\begin{align*}
g_t=&-\kappa\varepsilon(1+t)^{-\kappa-1}h^d
+\varepsilon(1+t)^{-\kappa}dh^{d-1}\frac{\beta|x|^2}{(1+t)^{\beta+1}},\\
\nabla g^m=&-\varepsilon^m(1+t)^{-m\kappa}mdh^{md-1}\frac{2x}{(1+t)^{\beta}},\\
\Delta g^m=&\varepsilon^m(1+t)^{-m\kappa}md(md-1)h^{md-2}%
\frac{4|x|^2}{(1+t)^{2\beta}}\\
&\quad-\varepsilon^m(1+t)^{-m\kappa}mdh^{md-1}%
\frac{2n}{(1+t)^{\beta}},
\end{align*}
for all $x\in A(t)$ and $t>0$.
According to the definition of $g$, we see that
$\pd{g}{\nu}\le0$ and $\pd{g^m}{\nu}\le0$ on $\partial\Omega$
since $\Omega$ is convex, and
$$g(x,0)=\varepsilon[(\eta-|x|^2)_+]^d\le \varepsilon_11_{B_{r}(x_0)}\le u_0(x),
\quad x\in\Omega,$$
provided that
\begin{equation} \label{eq-zcondi1}
\eta\le r^2, \quad \varepsilon\eta^d\le\varepsilon_1.
\end{equation}
In order to find a weak lower solution $g$, we only need to check the following
differential inequality on $A(t)$
\begin{align} \label{eq-zweak}
\pd{g}{t}\le\Delta g^m-\nabla\cdot(g^m\phi(g)\nabla v)
+\mu g^\delta(1-g), \quad x\in A(t), ~t>0.
\end{align}
Since $g(x,t)\le\varepsilon\eta^d\le\varepsilon_1\le1/2$, we see that
$\mu g^\delta(1-g)\ge\mu g^\delta/2$ for all $x\in\Omega$ and $t\ge0$.
Further,
\begin{align*}
|\nabla\cdot(g^m\phi(g)\nabla v)|
&\le g^m|\phi(g)||\Delta v|+|mg^{m-1}\phi(g)+g^m\phi'(g)||\nabla g||\nabla v|\\
&\le g^m\|\Delta v\|_{L^\infty(\Omega\times\mathbb R^+)}
+(m+1)|\nabla g^m|\cdot\|\nabla v\|_{L^\infty(\Omega\times\mathbb R^+)}.
\end{align*}
We denote $C_1=\|\nabla v\|_{L^\infty(\Omega\times\mathbb R^+)}$
and $C_2=\|\Delta v\|_{L^\infty(\Omega\times\mathbb R^+)}$ for convenience,
since they are bounded according to Theorem \ref{th-regular}.
A sufficient condition of inequality \eqref{eq-zweak} is
\begin{align} \nonumber
\varepsilon(1+t)^{-\kappa}&dh^{d-1}\frac{\beta|x|^2}{(1+t)^{\beta+1}}
+\varepsilon^m(1+t)^{-m\kappa}mdh^{md-1}%
\frac{2n}{(1+t)^{\beta}}\\ \nonumber
&\qquad+C_2\varepsilon^m(1+t)^{-m\kappa}h^{md}
+(m+1)C_1\varepsilon^m(1+t)^{-m\kappa}mdh^{md-1}\frac{2|x|}{(1+t)^{\beta}}
\\ \nonumber
&\le
\kappa\varepsilon(1+t)^{-\kappa-1}h^d
+\varepsilon^m(1+t)^{-m\kappa}md(md-1)h^{md-2}%
\frac{4|x|^2}{(1+t)^{2\beta}}
\\ \label{eq-zlower}
&\qquad+\frac{\mu}{2} \varepsilon^\delta(1+t)^{-\kappa\delta}h^{d\delta},
\quad x\in A(t), ~t>0.
\end{align}
As we have chosen $d=1/(m-1)$ and $\kappa=(1-\beta)/(m-1)$, we rewrite
\eqref{eq-zlower} into
\begin{align} \nonumber
\frac{\varepsilon\beta}{m-1}&\frac{|x|^2}{(1+t)^{\beta}}
+2n\frac{m}{m-1}\varepsilon^m h\\ \nonumber
&\qquad+C_2\varepsilon^m(1+t)^{\beta}h^2
+2(m+1)C_1\varepsilon^m\frac{m}{m-1}h|x|
\\ \nonumber
&\le
\kappa\varepsilon h
+\varepsilon^m\frac{m}{(m-1)^2}\frac{4|x|^2}{(1+t)^{\beta}}
\\
&\qquad+\frac{\mu}{2} \varepsilon^\delta(1+t)^{-\kappa\delta+\kappa+1}%
h^{d\delta-d+1},
\quad x\in A(t), ~t>0.
\end{align}
Let $\varepsilon$, $\beta$ and $\eta$ be chosen such that
\begin{align} \label{eq-zcondi2}
\begin{cases}
\varepsilon\beta\le 4\varepsilon^m\frac{m}{m-1}, \\
2n\frac{m}{m-1}\varepsilon^m\le\frac{1}{2}\kappa\varepsilon,\\
2mC_1\varepsilon^m|x|
\le\frac{1}{2}\kappa\varepsilon, \\
C_2\varepsilon^mh^{d+1-d\delta}
\le\frac{\mu}{2} \varepsilon^\delta(1+t)^{-\kappa\delta+\kappa+1-\beta},
\quad x\in A(t), ~t>0.
\end{cases}
\end{align}
Since $1\le\delta<m$, $\beta\in(0,1/2)$, $\kappa=(1-\beta)/(m-1)\ge1/[2(m-1)]$,
$h\le1/2$ and $|x|\le\text{diam}\Omega$,
we see that $d+1-d\delta=d(m-\delta)>0$,
$-\kappa\delta+\kappa+1-\beta=(m-\delta)\kappa>0$.
Thus, for \eqref{eq-zcondi1} and \eqref{eq-zcondi2},
it suffices to choose $\eta=r^2$,
$$\varepsilon=
\min\Big\{\Big(\frac{1}{8nm}\Big)^\frac{1}{m-1},
\Big(\frac{1}{8m(m-1)C_1\text{diam}\Omega}\Big)^\frac{1}{m-1},
\frac{\varepsilon_1}{r^{2d}},
\Big(\frac{\mu}{2C_2}\Big)^\frac{1}{m-\delta}\Big\},$$
and then $\beta=4\varepsilon^{m-1}m/(m-1)$.

Now, we find a weak lower solution with expanding support and
comparison principle Lemma \ref{le-comparisonprinciple} implies
$$u(x,t)\ge g(x,t)=\varepsilon(1+t)^{-\kappa}
\Big[\Big(\eta-\frac{|x-x_0|^2}{(1+t)^\beta}\Big)_+\Big]^d,
\quad x\in\Omega,~t>0.$$
There exists a $t_0$ such that
$$\eta-\frac{|x-x_0|^2}{(1+t_0)^\beta}\ge \frac{\eta}{2},
\quad x\in\Omega,$$
and thus
$$u(x,t_0)\ge g(x,t_0)\ge\varepsilon(1+t_0)^{-\kappa}
\Big(\frac{\eta}{2}\Big)^d,
\quad x\in\Omega.$$
Next, we construct another constant lower solution
$$\underline u(x,t)\equiv\varepsilon_0, \quad x\in\Omega, ~t>t_0,$$
with $0<\varepsilon_0\le\varepsilon(1+t_0)^{-\kappa}({\eta}/{2})^d\le1/2$
to be determined.
Clearly, $\pd{\underline u}{\nu}=0$ on $\partial\Omega$.
We only need to check the following differential inequality
$$0\le-\varepsilon_0^m\phi(\varepsilon_0)\Delta v(x,t)
+\mu \varepsilon_0^\delta(1-\varepsilon_0), \quad x\in\Omega, ~t>t_0,$$
which is valid if we further let
$$\varepsilon_0\le
\Big(\frac{\mu}{2\|\Delta v\|_{L^\infty(\Omega\times(t_0,+\infty))}}
\Big)^\frac{1}{m-\delta},$$
since $\delta<m$ and $\Delta v$ is uniformly bounded
according to Theorem \ref{th-regular}.
Applying the comparison principle Lemma \ref{le-comparisonprinciple} again,
we find
$$u(x,t)\ge\underline u(x,t)\equiv\varepsilon_0,
\quad x\in\Omega,~t>t_0.$$
This completes the proof.
$\hfill\Box$

{\it\bfseries Remark.}
It is interesting to compare the self similar weak lower solution $g(x,t)$
in the proof of Lemma \ref{le-lower} to the Barenblatt solution of porous medium equation
$$B(x,t)=(1+t)^{-k}\Big[\Big(1-\frac{k(m-1)}{2mn}\frac{|x|^2}{(1+t)^{2k/n}}
\Big)_+\Big]^\frac{1}{m-1},$$
with $k=1/(m-1+2/n)$.
The Barenblatt solution $B(x,t)$ is decaying at the rate $(1+t)^{-1/(m-1+2/n)}$
in $L^\infty(\mathbb R^n)$ and the support is expanding at the rate $(1+t)^{2k/n}$.
While the self similar weak lower solution $g(x,t)$
is decaying at the rate $(1+t)^{-(1-\beta)/(m-1)}$ and its support is
expanding at the rate $(1+t)^\beta$.
Here in the proof we have selected $\beta>0$ sufficiently small,
which means the support of $g$ is expanding with a much slower rate
and the maximum of $g$ is decaying at a slightly faster rate.

\noindent
{\it\bfseries Proof of Theorem \ref{th-lower}}\quad
This has been proved in Lemma \ref{le-lower}.
$\hfill\Box$

After proving the support expanding property of
the first equation in \eqref{eq-model},
which is a degenerate diffusion equation,
we can deduce the following convergence properties of all components.

\begin{lemma} \label{le-converge}
Let $(u,v,w,z)$ be a globally bounded weak solution of \eqref{eq-model}
with the first component initial data $u_0\ge0$, $u_0\not\equiv0$
and $1\le\delta<m$.
Then there exist constants $C_1,C_2>0$
and $c_1,c_2>0$ independent of $t$ such that
\begin{align*}
\|w(\cdot,t)\|_{L^\infty(\Omega)}
+\|\nabla w(\cdot,t)\|_{L^\infty(\Omega)}\le C_1e^{-c_1t},
\end{align*}
and
\begin{align*}
\|v(\cdot,t)-(\overline v_0+\overline w_0)\|_{L^\infty(\Omega)}
+\|\nabla v(\cdot,t)\|_{L^\infty(\Omega)}
+\|\Delta v(\cdot,t)\|_{L^\infty(\Omega)}\le C_2e^{-c_2t},
\end{align*}
for all $t>0$, where $\overline f=\int_\Omega fdx/|\Omega|$.
\end{lemma}
{\it\bfseries Proof.}
Applying Lemma \ref{le-intv}, we see that
\begin{align*}
\int_0^tz(x,t)ds&\ge C\int_0^t\int_\Omega u(y,s)dyds-C\\
&\ge C\int_{t_0}^t\int_\Omega u(y,s)dyds-C\\
&\ge C|\Omega|\varepsilon_0(t-t_0)-C\\
&\ge c_1t-C, \quad x\in\Omega,~t>t_0,
\end{align*}
since $u(x,t)\ge\varepsilon_0$ for $x\in\Omega$ and $t>t_0$
according to Lemma \ref{le-lower}.
Therefore,
\begin{align} \label{eq-zwdecay}
w(x,t)=w_0(x)e^{-\int_0^tz(x,s)ds}
\le w_0(x)e^{-c_1t+C}\le C_1e^{-c_1t}, \quad x\in\Omega,~t>t_0.
\end{align}
This is also valid for $t\in(0,t_0)$ upon enlarging $C_1$ if necessary
and hereafter we only need to prove this lemma for $t>t_0$.
We also have
\begin{align*}
|\nabla w(x,t)|&=|\nabla w_0(x)|e^{-\int_0^tz(x,s)ds}
+w_0(x)e^{-\int_0^tz(x,s)ds}\Big|\int_0^t\nabla z(x,s)ds\Big|\\
&\le Ce^{-c_1t}+Ce^{-c_1t}t
\le C_1'e^{-c_1't}, \quad x\in\Omega,~t>t_0,
\end{align*}
with $0<c_1'<c_1$.
We may write $C_1'$ and $c_1'$ as $C_1$ and $c_1$ for simplicity.
Therefore,
\begin{align*}
|\nabla (wz)(x,t)|\le|z\nabla w(x,t)|+|w\nabla z(x,t)|
\le Ce^{-c_1t}, \quad x\in\Omega,~t>t_0,
\end{align*}
It follows form the second equation in \eqref{eq-model} that
$$v(x,t)=e^{t\Delta}v_0+\int_0^te^{(t-s)\Delta}(wz)(\cdot,s)ds,
\quad t>0,$$
and
$$\nabla v(x,t)=e^{t\Delta}\nabla v_0+\int_0^te^{(t-s)\Delta}\nabla(wz)(\cdot,s)ds,
\quad t>0,$$
Using the standard $L^p-L^q$ type estimate for $\Delta v$, we get
\begin{align*}
\|\Delta v(x,t)\|_{L^\infty(\Omega)}
\le& \|\nabla e^{t\Delta}|\nabla v_0|\|_{L^\infty(\Omega)}
+\int_0^t\|\nabla e^{(t-s)\Delta}|\nabla(wz)(x,s)|\|_{L^\infty(\Omega)}ds\\
\le& C(1+t^{-\frac{1}{2}})e^{-\lambda_1t}\|\nabla v_0\|_{L^\infty(\Omega)}\\
&\qquad+C\int_0^t(1+(t-s)^{-\frac{1}{2}})e^{-\lambda_1(t-s)}
\|\nabla(wz)(\cdot,s)\|_{L^\infty(\Omega)}\\
\le& Ce^{-\lambda_1t}
+C\int_0^t(1+(t-s)^{-\frac{1}{2}})e^{-\lambda_1(t-s)}e^{-c_1s}ds\\
\le& C_2e^{-c_2t}, \quad x\in\Omega,~t>t_0,
\end{align*}
where $\lambda_1>0$ is the first nonzero eigenvalue of $-\Delta$ with homogeneous
Neumann boundary condition.
The $L^\infty$ estimate of $\nabla v$ can be deduced in a similar way.
In the proof of Lemma \ref{le-wv}, we have obtained
$$\int_\Omega (v(x,t)+w(x,t))dx\equiv \int_\Omega (v_0(x)+w_0(x))dx,$$
which is the same as the estimate of $v_\varepsilon+w_\varepsilon$.
It follows from \eqref{eq-zwdecay} that
$w(x,t)$ is decaying to zero exponentially.
This implies that
$$\overline v(t)=\frac{1}{|\Omega|}\int_\Omega v(x,t)dx$$
is converging to $\overline v_0+\overline w_0$ exponentially.
A Poincar\'e type inequality shows
$$\|v(x,t)-\overline v(t)\|_{L^\infty(\Omega)}
\le C\|\nabla v(x,t)\|_{L^\infty(\Omega)}\le Ce^{-c_2t}.$$
Therefore,
\begin{align*}
\|v(x,t)-(\overline v_0+\overline w_0)\|_{L^\infty(\Omega)}
&\le \|v(x,t)-\overline v(t)\|_{L^\infty(\Omega)}
+\|\overline v(t)-(\overline v_0+\overline w_0)\|_{L^\infty(\Omega)}\\
&\le \|v(x,t)-\overline v(t)\|_{L^\infty(\Omega)}
+\|\overline w(t)\|_{L^\infty(\Omega)}\\
&\le Ce^{-c_2't},
\quad x\in\Omega,~t>t_0,
\end{align*}
The proof is completed.
$\hfill\Box$

\begin{lemma} \label{le-ode}
For constants $C,c>0$ and $m>1$,
the local solution $g$ of the following ODE
\begin{equation*}
\begin{cases}
g'(t)=Ce^{-ct}g^m, \quad t>0,\\
g(0)=g_0>0,
\end{cases}
\end{equation*}
blows up in finite time if $c/C<(m-1)g_0^{m-1}$,
while remains bounded if $c/C>(m-1)g_0^{m-1}$.
\end{lemma}
{\it\bfseries Proof.}
There holds
$$\frac{-1}{m-1}\Big(\frac{1}{g^{m-1}}\Big)'=Ce^{-ct}, \quad t>0.$$
Integrating over $(0,t)$ shows
$$\frac{1}{m-1}\Big(\frac{1}{g_0^{m-1}}-\frac{1}{g^{m-1}(t)}\Big)
=\frac{C}{c}(1-e^{-ct}).$$
A simple analysis completes this proof.
$\hfill\Box$

\begin{lemma} \label{le-converge2}
Let $(u,v,w,z)$ be a globally bounded weak solution of \eqref{eq-model}
with the first component initial data $u_0\ge0$, $u_0\not\equiv0$
and $1\le\delta<m$.
Then there exist constants $C_3>0$
and $c_3>0$ independent of $t$ such that
\begin{align*}
\|u(\cdot,t)-1\|_{L^\infty(\Omega)}\le C_3e^{-c_3t},
\end{align*}
for all $t>0$.
\end{lemma}
{\it\bfseries Proof.}
Lemma \ref{le-lower} implies that $u(x,t)\ge\varepsilon_0$
for $x\in\Omega$ and $t>t_0$.
It suffices to prove this lemma for $t\ge t_1$
with some fixed $t_1\ge t_0$ to be determined.
We use upper and lower solution method to achieve this.
Let $u_1(t)$ and $u_2(t)$ be
one pair of the solutions of the following ODE
\begin{align} \label{eq-zupper}
\begin{cases}
u_1'(t)\ge u_1^m\|\Delta v(\cdot,t)\|_{L^\infty(\Omega)}
+\mu u_1^\delta(1-u_1), \\
u_2'(t)\le -u_2^m\|\Delta v(\cdot,t)\|_{L^\infty(\Omega)}
+\mu u_2^\delta(1-u_2), \quad t>t_1, \\
u_1(t_1)\ge \|u(\cdot,t_1)\|_{L^\infty(\Omega)}, \\
u_2(t_1)\le \varepsilon_0.
\end{cases}
\end{align}
Lemma \ref{le-comparisonprinciple} shows that
$$u_1(t)\ge u(x,t)\ge u_2(t), \quad x\in\Omega,~t>t_0.$$
We only need to find one pair of $(u_1,u_2)$ such that
$u_1$ and $u_2$ both converge to $1$ exponentially.
A sufficient condition of \eqref{eq-zupper} is
\begin{align} \label{eq-zupper2}
\begin{cases}
u_1'(t)=C_2e^{-c_2t}u_1^m+\mu u_1^\delta(1-u_1), \\
u_2'(t)=-C_2e^{-c_2t}u_2^m+\mu u_2^\delta(1-u_2), \quad t>t_1, \\
u_1(t_1)=\|u(\cdot,t_1)\|_{L^\infty(\Omega)}+1, \\
u_2(t_1)=\varepsilon_0,
\end{cases}
\end{align}
since $\|\Delta v(\cdot,t)\|_{L^\infty(\Omega)}\le C_2e^{-c_2t}$
according to Lemma \ref{le-converge}.
We note that we can choose $t_1$ sufficiently large
such that
$$\frac{c_2}{C_2e^{-c_2t_1}}>2(m-1)
\Big(\sup_{t>0}\|u(\cdot,t)\|_{L^\infty(\Omega)}+1\Big)^{m-1}.$$
Lemma \ref{le-ode} implies that $u_1(t)$ is uniformly bounded by some constant $C$.
And a simple ODE comparison shows that $u_1(t)>1$ for all $t>t_1$.
Therefore,
\begin{align*}
\begin{cases}
u_1'(t)\le C^m C_2e^{-c_2t}+\mu \varepsilon_0^\delta(1-u_1), \quad t>t_1, \\
u_1(t_1)=\|u(\cdot,t_1)\|_{L^\infty(\Omega)}+1.
\end{cases}
\end{align*}
We see that $u_1(t)$ is an upper solution of $u(x,t)$
and an upper solution of $u_1(t)$ is $\overline u_1(t)$ such that
\begin{align} \label{eq-zu1upper}
\begin{cases}
\overline u_1'(t)=C^m C_2e^{-c_2t}
+\mu \varepsilon_0^\delta(1-\overline u_1), \quad t>t_1, \\
\overline u_1(t_1)=\|u(\cdot,t_1)\|_{L^\infty(\Omega)}+1,
\end{cases}
\end{align}
which can be solved as
\begin{align*}
\overline u_1(t)
&=1+e^{-\mu\varepsilon_0^\delta (t-t_1)}(\|u(\cdot,t_1)\|_{L^\infty(\Omega)}+1)
+C^m C_2\int_{t_1}^te^{-\mu\varepsilon_0^\delta (t-s)}e^{-c_2s}ds
-e^{-\mu\varepsilon_0^\delta (t-t_1)}\\
&\le 1+e^{-\mu\varepsilon_0^\delta (t-t_1)}\|u(\cdot,t_1)\|_{L^\infty(\Omega)}
+C^m C_2Ce^{-\min\{\mu\varepsilon_0^\delta,c_2\}t/2}, \quad t>t_1.
\end{align*}
On the other hand, the lower solution of $u(x,t)$ satisfies
\begin{align*}
\begin{cases}
u_2'(t)=-C_2e^{-c_2t}u_2^m+\mu u_2^\delta(1-u_2), \quad t>t_1, \\
u_2(t_1)=\varepsilon_0.
\end{cases}
\end{align*}
We note that we can choose $t_1$ sufficiently large that
$$C_2e^{-c_2t}\varepsilon_0^m\le\mu \varepsilon_0^\delta(1-\varepsilon_0).$$
An ODE comparison shows that $\varepsilon_0\le u_2(t)<1$ for all $t>t_1$ and
\begin{align*}
\begin{cases}
u_2'(t)\ge-C_2e^{-c_2t}+\mu \varepsilon_0^\delta(1-u_2), \quad t>t_1, \\
u_2(t_1)=\varepsilon_0.
\end{cases}
\end{align*}
We see that $u_2(t)$ is a lower solution of $u(x,t)$
and a lower solution of $u_2(t)$ is $\underline u_2(t)$ such that
\begin{align*}
\begin{cases}
\underline u_2'(t)=-C_2e^{-c_2t}
+\mu \varepsilon_0^\delta(1-\underline u_2), \quad t>t_1, \\
\underline u_2(t_1)=\varepsilon_0.
\end{cases}
\end{align*}
This can also be solved as
\begin{align*}
\underline u_2(t)
&=1+e^{-\mu\varepsilon_0^\delta (t-t_1)}\varepsilon_0
-C_2\int_{t_1}^te^{-\mu\varepsilon_0^\delta (t-s)}e^{-c_2s}ds
-e^{-\mu\varepsilon_0^\delta (t-t_1)}\\
&\ge 1-e^{-\mu\varepsilon_0^\delta (t-t_1)}
-C_2Ce^{-\min\{\mu\varepsilon_1^\delta,c_2\}t/2}, \quad t>t_1.
\end{align*}
Thus, we conclude
$$\underline u_2(t)\le u_t(t)\le u(x,t)\le u_t(t)\le\overline u_1(t),
\quad t>t_1,$$
and $\underline u_2(t)$, $\overline u_1(t)$ converge to $1$ exponentially.
$\hfill\Box$

\begin{lemma} \label{le-converge3}
Let $(u,v,w,z)$ be a globally bounded weak solution of \eqref{eq-model}
with the first component initial data $u_0\ge0$, $u_0\not\equiv0$
and $1\le\delta<m$.
Then there exist constants $C_4>0$
and $c_4>0$ independent of $t$ such that
\begin{align*}
\|z(\cdot,t)-1\|_{L^\infty(\Omega)}\le C_4e^{-c_4t},
\end{align*}
for all $t>0$.
\end{lemma}
{\it\bfseries Proof.}
From the fourth equation in \eqref{eq-model}, we have
$$z(x,t)=e^{t(\Delta-1)}z_0
+\int_0^te^{(t-s)(\Delta-1))}u(\cdot,s)ds,
\quad t>0.$$
We note that
$$\int_0^te^{(t-s)(\Delta-1))}1ds=1-e^{-t},$$
which can be deduced by solving the ODE $z'=-z+1$ with $z(0)=0$.
Therefore,
\begin{align*}
\|z(x,t)-1\|_{L^\infty(\Omega)}
&\le \|e^{t(\Delta-1)}z_0\|_{L^\infty(\Omega)}+
\int_0^t\|e^{(t-s)(\Delta-1))}(u(\cdot,s)-1)\|_{L^\infty(\Omega)}ds
+e^{-t}\\
&\le Ce^{-t}(\|z_0\|_{L^\infty(\Omega)}+1)
+C\int_0^te^{-(t-s)}\|(u(\cdot,s)-1)\|_{L^\infty(\Omega)}ds\\
&\le Ce^{-t}+CC_3\int_0^te^{-(t-s)}e^{-c_3s}ds\\
&\le C_4e^{-c_4t}, \quad t>0.
\end{align*}
The proof is completed.
$\hfill\Box$

\noindent
{\it\bfseries Proof of Theorem \ref{th-asymp}}\quad
This is proved by collecting Lemma \ref{le-lower},
Lemma \ref{le-converge}, Lemma \ref{le-converge2}
and Lemma \ref{le-converge3}.
$\hfill\Box$

Finally, we construct a self similar upper solution with expanding support
to prove Theorem \ref{th-upper}.
We note that for constructing a weak upper solution for the heat equation,
one should replace the cut-off composite function $(\cdot)_+$ by $(\cdot)_-$.
But here for the degenerate porous medium type equation
and the self similar function of the form
$g=[(1-|x|^2)_+]^d$ with $md>1$,
we can check that $\nabla g^m$ is continuous and
$\Delta g^m\in L^q(\Omega)$ for some $q>1$.
This shows that the differential inequality for an upper solution
only need to be valid almost everywhere,
without the possible Radon measures on the boundary of its support.

\begin{lemma} \label{le-upper}
Let $(u,v,w,z)$ be a globally bounded weak solution of \eqref{eq-model}.
We further assume that
$$\text{supp }u_0\subset\overline B_{r_0}(x_0)\subset\Omega,$$
for some $r_0>0$ and $x_0\in\Omega$.
Define a function
$$g(x,t)=\varepsilon(\tau+t)^\sigma
\Big[\Big(\eta-\frac{|x-x_0|^2}{(\tau+t)^\beta}\Big)_+\Big]^d,
\quad x\in\Omega,~t\ge0,$$
where $d=1/(m-1)$, $\beta>0$, $\sigma>0$,
$\varepsilon>0$, $\eta>0$, $\tau\in(0,1)$.
Then by appropriately selecting $\beta$, $\sigma$
$\varepsilon$, $\eta$ and $\tau$,
the support of $g(x,t)$ is contained in $\Omega$ for
$t\in(0,t_0)$ with some $t_0>0$ and
the function $g(x,t)$ can be
an upper solution of the first equation in \eqref{eq-model}
on $\Omega\times(0,t_0)$, that is,
\begin{align*}
\begin{cases}
\displaystyle
\pd{g}{t}\ge\Delta g^m-\nabla\cdot(g^m\phi(g)\nabla v)
+\mu g^\delta(1-g),\qquad &x\in\Omega, t\in(0,t_0), \\[3mm]
\displaystyle
\pd{g}{\nu}\ge0, \qquad &x\in\partial\Omega, t\in(0,t_0),\\[2mm]
g(x,0)\ge u_0(x)\ge0, \qquad &x\in\Omega,
\end{cases}
\end{align*}
in the sense that the following inequality
\begin{align*}
\int_0^{t_0}\int_\Omega g\varphi_tdxdt&
+\int_\Omega g(x,0)\varphi(x,0) dx
\le\int_0^{t_0}\int_\Omega \nabla g^m\cdot\nabla\varphi dxdt\\
&\qquad-\int_0^{t_0}\int_\Omega g^m\phi(g)\nabla v\cdot\nabla\varphi dxdt
-\int_0^{t_0}\int_\Omega\mu g^\delta(1-g)\varphi dxdt,
\end{align*}
holds for all test functions
$0\le\varphi\in L^2((0,t_0);W^{1,2}(\Omega))\cap W^{1,2}((0,t_0);L^2(\Omega))$ with
$\varphi(x,t_0)=0$ on $\Omega$
and $g(x,0)\ge u_0(x)\ge0$ on $\Omega$.
Therefore, $u(x,t)\le g(x,t)$ and there exist
a family of monotone increasing open sets
$\{A(t)\}_{t\in(0,t_0)}$ such that
$$\text{supp }u(\cdot,t)\subset\overline A(t)\subset\Omega, \quad t\in(0,t_0),$$
and $\partial A(t)$ has a finite derivative with respect to $t$.
\end{lemma}
{\it\bfseries Proof.}
For simplicity, we let
$$h(x,t)=\Big(\eta-\frac{|x-x_0|^2}{(\tau+t)^\beta}\Big)_+, \quad x\in\Omega, ~t\ge0,$$
and
$$A(t)=\Big\{x\in\Omega;\frac{|x-x_0|^2}{(\tau+t)^\beta}<\eta\Big\},
\quad t\ge0.$$
Since $u_0\in C(\overline\Omega)$ and
$\text{supp }u_0\subset\overline B_{r_0}(x_0)\subset\Omega$,
we see that there exist $r_1>r_0$ and $\varepsilon_1>0$ such that
$B_{r_1}(x_0)\subset\subset\Omega$ and
$u_0(x)\le\varepsilon_1$ for all $x\in\Omega$.
Without loss of generality, we may assume that $x_0=0$.
Straightforward computation shows that
\begin{align*}
g_t=&\sigma\varepsilon(\tau+t)^{\sigma-1}h^d
+\varepsilon(\tau+t)^{\sigma}dh^{d-1}\frac{\beta|x|^2}{(\tau+t)^{\beta+1}},\\
\nabla g^m=&-\varepsilon^m(\tau+t)^{m\sigma}mdh^{md-1}\frac{2x}{(\tau+t)^{\beta}},\\
\Delta g^m=&\varepsilon^m(\tau+t)^{m\sigma}md(md-1)h^{md-2}%
\frac{4|x|^2}{(\tau+t)^{2\beta}}-\varepsilon^m(\tau+t)^{m\sigma}mdh^{md-1}%
\frac{2n}{(\tau+t)^{\beta}},
\end{align*}
for all $x\in A(t)$ and $t>0$.
Let $\tau\in(0,1)$ to be determined and
\begin{equation} \label{eq-zuppercondi}
r_2=\frac{r_0+r_1}{2}, \quad \eta=\frac{r_2^2}{\tau^\beta},
\quad t_0=\min\Big\{\tau,\tau\Big(\Big(\frac{r_1}{r_2}\Big)^\frac{2}{\beta}-1
\Big)\Big\}.
\end{equation}
According to the definition of $g$, we see that
$A(0)=B_{r_2}(x_0)$, $\text{supp }u_0\subset\subset\overline  A(0)\subset\Omega$,
and $A(t_0)\subset B_{r_1}(x_0)\subset\subset\Omega$.
Therefore, $\pd{g}{\nu}=0$ and $\pd{g^m}{\nu}=0$ on $\partial\Omega$
for all $t\in(0,t_0)$, and
\begin{equation*}
g(x,0)=\varepsilon\tau^\sigma
\Big[\Big(\eta-\frac{|x-x_0|^2}{\tau^\beta}\Big)_+\Big]^d
\ge\varepsilon\tau^\sigma\Big(\frac{r_2^2}{\tau^\beta}-\frac{r_0^2}{\tau^\beta}\Big)^d
\cdot1_{B_{r_0}(x_0)}
\ge\varepsilon_11_{B_{r_0}(x_0)}\ge u_0(x),
\quad x\in\Omega,
\end{equation*}
provided that
\begin{equation} \label{eq-zuppercondi1}
\varepsilon\tau^\sigma\Big(\frac{r_2^2}{\tau^\beta}-\frac{r_0^2}{\tau^\beta}\Big)^d
\ge\varepsilon_1.
\end{equation}
In order to find a weak lower solution $g$, we only need to check the following
differential inequality on $A(t)$
\begin{align} \label{eq-zupperweak}
\pd{g}{t}\ge\Delta g^m-\nabla\cdot(g^m\phi(g)\nabla v)
+\mu g^\delta(1-g), \quad x\in A(t), ~t\in(0,t_0).
\end{align}
Since $0\le g\le\varepsilon\eta^d$,
we see that $\mu g^\delta(1-g)\le\mu g^\delta$ for all $x\in\Omega$ and $t\ge0$.
Further,
\begin{align*}
|\nabla\cdot(g^m\phi(g)\nabla v)|
&\le g^m|\phi(g)||\Delta v|+|mg^{m-1}\phi(g)+g^m\phi'(g)||\nabla g||\nabla v|\\
&\le g^m\|\Delta v\|_{L^\infty(\Omega\times\mathbb R^+)}
+(m+\varepsilon\tau^\sigma\eta^d)|\nabla g^m|
\cdot\|\nabla v\|_{L^\infty(\Omega\times\mathbb R^+)}.
\end{align*}
We denote $C_1=\|\nabla v\|_{L^\infty(\Omega\times\mathbb R^+)}$
and $C_2=\|\Delta v\|_{L^\infty(\Omega\times\mathbb R^+)}$ for convenience,
since they are bounded according to Theorem \ref{th-regular}.
A sufficient condition of inequality \eqref{eq-zupperweak} is
\begin{align} \nonumber
\sigma\varepsilon&(\tau+t)^{\sigma-1}h^d
+\varepsilon(\tau+t)^{\sigma}dh^{d-1}\frac{\beta|x|^2}{(\tau+t)^{\beta+1}}
+\varepsilon^m(\tau+t)^{m\sigma}mdh^{md-1}\frac{2n}{(\tau+t)^{\beta}}
\\ \nonumber
\ge&C_2\varepsilon^m(\tau+t)^{m\sigma}h^{md}
+(m+\varepsilon\tau^\sigma\eta^d)
C_1\varepsilon^m(\tau+t)^{m\sigma}mdh^{md-1}\frac{2|x|}{(\tau+t)^{\beta}}
\\ \label{eq-zuppersimilar}
&+\varepsilon^m(\tau+t)^{m\sigma}md(md-1)h^{md-2}%
\frac{4|x|^2}{(\tau+t)^{2\beta}}
+\mu\varepsilon^\delta (\tau+t)^{\delta\sigma}h^{d\delta},
\quad x\in A(t), ~t\in(0,t_0).
\end{align}
As we have chosen $d=1/(m-1)$, we rewrite \eqref{eq-zuppersimilar} into
\begin{align*}
\sigma\varepsilon&(\tau+t)^{\sigma-1}h
+\frac{\varepsilon\beta}{m-1}(\tau+t)^{\sigma}\frac{|x|^2}{(\tau+t)^{\beta+1}}
+2n\frac{m}{m-1}\varepsilon^m(\tau+t)^{m\sigma}\frac{h}{(\tau+t)^{\beta}}
\\
\ge&C_2\varepsilon^m(\tau+t)^{m\sigma}h^2
+2(m+\varepsilon\tau^\sigma\eta^d)
C_1\varepsilon^m(\tau+t)^{m\sigma}mdh\frac{|x|}{(\tau+t)^{\beta}}
\\
&+\frac{m}{(m-1)^2}\varepsilon^m(\tau+t)^{m\sigma}%
\frac{4|x|^2}{(\tau+t)^{2\beta}}
+\mu\varepsilon^\delta (\tau+t)^{\delta\sigma}h^{d\delta-d+1},
\quad x\in A(t), ~t\in(0,t_0).
\end{align*}
Let $\varepsilon$, $\beta$, $\sigma$ and $\tau$ be chosen such that
\begin{align} \label{eq-zuppercondi2}
\begin{cases}
\displaystyle
\frac{1}{2}\frac{\varepsilon\beta}{m-1}
(\tau+t)^{\sigma}\frac{|x|^2}{(\tau+t)^{\beta+1}}
\ge\frac{m}{(m-1)^2}\varepsilon^m(\tau+t)^{m\sigma}%
\frac{4|x|^2}{(\tau+t)^{2\beta}}, \\
\displaystyle
\frac{1}{3}\sigma\varepsilon(\tau+t)^{\sigma-1}h
\ge C_2\varepsilon^m(\tau+t)^{m\sigma}h^2,\\
\displaystyle
\frac{1}{3}\sigma\varepsilon(\tau+t)^{\sigma-1}h
\ge\mu\varepsilon^\delta (\tau+t)^{\delta\sigma}h^{d\delta-d+1}, \\
\displaystyle
\frac{1}{2}\frac{\varepsilon\beta}{m-1}
(\tau+t)^{\sigma}\frac{|x|^2}{(\tau+t)^{\beta+1}}
+\frac{1}{3}\sigma\varepsilon(\tau+t)^{\sigma-1}h
\ge2(m+\varepsilon\tau^\sigma\eta^d)
C_1\varepsilon^m(\tau+t)^{m\sigma}mdh\frac{|x|}{(\tau+t)^{\beta}},
\\[2mm]
\hfill\quad x\in A(t), ~t\in(0,t_0).
\end{cases}
\end{align}
We have the following estimate
\begin{align*}
&2(m+\varepsilon\tau^\sigma\eta^d)
C_1\varepsilon^m(\tau+t)^{m\sigma}mdh\frac{|x|}{(\tau+t)^{\beta}}\\
\le &\frac{m}{(m-1)^2}\varepsilon^m(\tau+t)^{m\sigma}%
\frac{4|x|^2}{(\tau+t)^{2\beta}}
+(m+\varepsilon\tau^\sigma\eta^d)^2
C_1^2m\varepsilon^m(\tau+t)^{m\sigma}h^2,
\quad x\in A(t), ~t\in(0,t_0).
\end{align*}
Therefore, a sufficient condition of \eqref{eq-zuppercondi2} is
\begin{align} \label{eq-zuppercondi3}
\begin{cases}
(m-1)\beta\ge8m\varepsilon^{m-1}(\tau+t)^{(m-1)\sigma-\beta+1}, \\
2\sigma/3
\ge (C_2+(m+\varepsilon\tau^\sigma\eta^d)^2C_1^2m)
\varepsilon^{m-1}(\tau+t)^{(m-1)\sigma+1}h,\\
\sigma/3
\ge\mu\varepsilon^{\delta-1}(\tau+t)^{(\delta-1)\sigma+1}h^{d(\delta-1)},
\quad x\in A(t), ~t\in(0,t_0).
\end{cases}
\end{align}
We note that $\eta$, $\tau$ and $t_0$ satisfy
the condition \eqref{eq-zuppercondi} and \eqref{eq-zuppercondi1},
and then $h\le\eta=r_2^2/\tau^\beta$,
$\tau+t\le\tau+t_0\le2\tau$,
$\varepsilon\tau^{\sigma-d\beta}(r_2^2-r_0^2)^d\ge\varepsilon_1$.
For $\tau\in(0,1)$, we choose
$$\varepsilon=\frac{\varepsilon_1}{\tau^{\sigma-d\beta}(r_2^2-r_0^2)^d}
:=C_3\tau^{d\beta-\sigma}.$$
Now, we only need to find $\tau\in(0,1)$ such that
\begin{align} \nonumber
\begin{cases}
(m-1)\beta\ge8mC_3^{m-1}2^{\max\{0,(m-1)\sigma-\beta+1\}}\tau, \\
2\sigma/3
\ge (C_2+(m+C_3r_2^{2d})^2C_1^2m)C_3^{m-1}2^{(m-1)\sigma+1}r_2^2\tau,\\
\sigma/3
\ge\mu C_3^{\delta-1}2^{(\delta-1)\sigma+1}r_2^{2d(\delta-1)}\tau.
\end{cases}
\end{align}
This can be done by selecting $\beta=1$, $\sigma=1$, and $\tau\in(0,1)$
sufficiently small.

The comparison principle Lemma \ref{le-comparisonprinciple} implies that
$u(x,t)\le g(x,t)$ for all $x\in\Omega$ and $t\in(0,t_0)$.
Thus,
$$\text{supp }u(\cdot,t)\subset\overline A(t)
=\{x\in\Omega;|x-x_0|^2<\eta(\tau+t)^\beta\},
\quad t\in(0,t_0),$$
and
$$\partial A(t)=\{x\in\Omega;|x-x_0|=\eta^\frac{1}{2}(\tau+t)^\frac{\beta}{2}\},
\quad t\in(0,t_0),$$
which has finite derivative with respect to $t$.
$\hfill\Box$

{\it\bfseries Remark.}
Similar to the weak lower solution in Lemma \ref{le-lower}, we
compare the self similar weak upper solution $g(x,t)$
in the proof of Lemma \ref{le-upper} to the Barenblatt solution of porous medium equation
$$B(x,t)=(1+t)^{-k}\Big[\Big(1-\frac{k(m-1)}{2mn}\frac{|x|^2}{(1+t)^{2k/n}}
\Big)_+\Big]^\frac{1}{m-1},$$
with $k=1/(m-1+2/n)$.
The Barenblatt solution $B(x,t)$ is decaying at the rate $(1+t)^{-1/(m-1+2/n)}$
in $L^\infty(\mathbb R^n)$ and the support is expanding at the rate $(1+t)^{2k/n}$.
As we have shown the support of the lower solution in Lemma \ref{le-lower}
is expanding with a much slower rate and decaying at a slightly faster rate.
Here, the upper solution
is increasing at the rate $(\tau+t)^{\sigma}$ and its support is
expanding at the rate $(\tau+t)^\beta$.
The increasing of $g(x,t)$ makes it possible to be an upper solution,
which can be seen from the proof.

From the proof of Lemma \ref{le-upper}, we can choose $\beta>0$
to as small as we want.
But we note that $\text{supp }u_0\subset\subset\text{supp }g(\cdot,0)$
and if we choose a smaller $\beta>0$,
then the parameters $\tau$ and $t_0$ are also smaller.
This shows if we let the upper solution expands slower,
then it may only be an upper solution for a smaller time interval.
Thus, the slower expanding upper solution $g(x,t)$
on a smaller time interval
does not contradict to the
possible feature that the solution $u(x,t)$ expands at a fixed rate
since $\text{supp }u_0\subset\subset\text{supp }g(\cdot,0)$ at the initial time.

\noindent
{\it\bfseries Proof of Theorem \ref{th-upper}}\quad
This has been proved in Lemma \ref{le-upper}.
$\hfill\Box$

\medskip

{\bf Acknowledgement}.  The research was initialed when M. Mei visited South China Normal University, and was finalized
when T. Xu and S. Ji visited McGill University. They would like to thank both universities for the great hospitality.
The research of S. Ji is
supported by the Fundamental Research Funds for the Central Universities (No.~2017BQ109),
, China Postdoctoral Science Foundation (No.~2017M610517)
and NSFC Grant No.~11701184.
The research of C. Jin was supported in part by NSFC grant No. 11471127,
Guangdong Natural Science Funds for Distinguished Young Scholar Grant No. 2015A030306029,
the Excellent Young Professors Program of Guangdong Province Grant No. HS2015007, and
Special Support Program of Guangdong Province of China.
The research of M. Mei was supported in part
by NSERC Grant RGPIN 354724-16, and FRQNT Grant No. 192571.
The research of J. Yin was supported in
part by NSFC Grant No. 11371153.


\end{document}